\documentclass[reqno]{amsart}
\usepackage{eucal}
\usepackage{a4wide}
\usepackage{graphicx} 
\usepackage[T1]{fontenc}      
\newtheorem{theo}{Theorem}
\newtheorem{prop}{Proposition}
\newtheorem{lm}{Lemma}
\newtheorem{lem}{Lemma}
\newtheorem{rem}{Remark}
\newtheorem{example}{Example}
\newtheorem{defi}{Definition}

\newtheorem{corr}{Corollary}
\newcommand{\mt}[1]{\ensuremath{\mathbf{#1}}}
\newcommand{\ms}[1]{\ensuremath{\mathbb{#1}}}
\newcommand{\lip}{Lipschitz }
\renewcommand{\div}{\mathrm{div}}
\newcommand{\Jac}{\mathrm{Jac}}
\newcommand{\Conv}{\mathrm{Conv}}
\begin{document}
\title[Hamiltonian approach to geodesic image matching]{Hamiltonian approach to geodesic image matching}
\author{Fran\c{c}ois-Xavier Vialard}
\address[Fran\c{c}ois-Xavier Vialard]{CMLA \\Ecole Normale
  Sup\'erieure de Cachan, CNRS, UniverSud\\ 61, avenue du Pr\'esident Wilson\\ F-94 235
  Cachan CEDEX}
\email{francois.xavier.vialard@normalesup.org}
\keywords{Variational calculus, energy minimization, Hamiltonian System, shape representation and recognition, geodesic, infinite dimensional riemannian manifolds, Lipschitz domain}
\subjclass{Primary: 58b10; Secondary: 49J45, 68T10}
\date{\today}
\begin{abstract}
This paper presents a generalization to image matching of the Hamiltonian approach
for planar curve matching developed in the context of group of diffeomorphisms. 
We propose an efficient framework to deal with discontinuous images
in any dimension, for example 2D or 3D. In this context, we give the structure of
the initial momentum (which happens to be decomposed in a smooth part and a singular part) thanks to a derivation lemma interesting in itself. 
The second part develops a Hamiltonian interpretation of the variational problem, derived from the optimal control theory point of view. 
\end{abstract}

\maketitle

\tableofcontents
\addtolength{\baselineskip}{0.3\baselineskip}

\section{Introduction}
This paper arose from the attempt to develop the multi-modal image matching in the framework 
of \emph{large deformation diffeomorphisms}.
Initiated by the work of Grenander, this context was deeply used since 
\cite{0855.57035}, especially with applications to computational anatomy. The method followed is the classical 
minization of an energy on the space of diffeomorphisms, which enables to compute geodesics on this space
and to derive the evolution equations. In most of the papers, the group of diffeomorphisms acts on the support
of the template; we add to this one a diffeomorphisms group action on the level set of the template. 
This action is a natural way to cope with the multi-modal matching and could be, in a certain way, compared to
the metamorphoses approach exposed in \cite{trouve05:_local_geomet_defor_templ}: the metamorphoses are another way
to act on the images but the goal is very different in our case.
Matching in our context, is to find a couple $(\eta,\phi)$ which minimizes the energy
\begin{equation} \label{minimisation}
E(\eta,\phi)=D(Id,(\eta,\phi))^2 + \frac{1}{\sigma^2} \| \eta \circ I_0 \circ \phi^{-1} -I_{targ}\|_{L^2}^2 ,
\end{equation}
with $I_0$ the initial function (or image) and $I_{targ}$ the target function, $\text{Id}$ is the identity map in the product 
of groups, $\sigma$ is a calibration parameter. 
The distance $D$ is obtained through a product of Riemannian metrics on the diffeomorphisms groups.
\\
All the complexity is then carried by the group of diffeomorphisms and its action: in the particular case
of landmark matching, the geodesics are well described. The problem is reduced in this case to understand the geodesic flow
on a finite dimensional Riemannian manifold. It should be also emphasized that this problem can be seen as an 
optimal control problem.
In \cite{BegIJCV}, numerical implementation of gradient based methods are strongly developed through a semi-Lagrangian
method for computing the geodesics.
A Hamiltonian formulation can be adopted 
to provide efficient applications and computations through the use of the conservation of momenta. 
In \cite{vmty04}, statistics are done on the initial momenta which is a relative 
signature of the target functions. The existence of geodesics from an initial momentum was deeply 
developed in \cite{trouve05:_local_geomet_defor_templ}, but this work dealt only with smooth functions for $I_0$ (essentially $H^1$)
however with a very large class of momenta. An attempt to understand the structure of the momentum for an initial discontinuous 
function was done in the matching of planar curves in \cite{hamcurves}. 
\\
We propose thereafter a framework to treat discontinuous
functions in any dimension: the main point is to derive the energy function in this context.
\\
Finally, we chose to give a Hamiltonian interpretation of the equations which is the proper way 
to handle the conservation of momentum. This formulation includes the work done in \cite{hamcurves} but
does not capture the landmark matching. The formulation we adopt gives a weak sense to the equations and
we prove existence and uniqueness for the weak Hamiltonian equations within a large set of initial data.
A word on the structure of initial data:
the article on planar matching (\cite{hamcurves}) focuses on Jordan curves. The main result is the existence for all time and
uniqueness of Hamiltonian flow. The initial data are roughly a Jordan curve for the position variable
and a vector field on this curve for the momentum. In our context, we choose four variables
$(I_0,\Sigma_0,p_0,P_0)$. $I_0$ is the initial function with a set of dicontinuities $\Sigma_0$,
$p_0$ is the momentum on the set $\Sigma_0$ and $P_0$ is the momentum for the smooth part of the initial function.
This is a natural way to understand the problem and 
the choice to keep the set of discontinuity as a position variable 
can lead to larger applications than only choose $I_0$ as position variable.

\vspace{0.3cm}

The paper is organized as follows. We start with a presentation of the framework underlying equation \eqref{minimisation}.
We present a key lemma concerning the data attachment term with respect to the $\eta$ and $\phi$ variables. Its proof
is postponed to the last section.
Then we derive the geodesics equations and ensure the existence of a solution for all time from an initial
momentum. In the second part of this paper, we give the weak formulation of the Hamiltonian equations, and deal with the
existence and uniqueness for this Hamiltonian formulation.

\section{Framework and notations}

\subsection{The space of discontinuous images}

\noindent
Let $n \geq 1$ and $M \subset \ms{R}^n$ a $C^1$ bounded open set diffeomorphic to the unit ball. 
\\
We denote by $BV(M)$ the set of functions of bounded variation. The reader is not supposed to have a broad knowledge 
of $BV$ functions. Below, we restrict ourselves to a subset of $BV$ functions which does not require the technical material of 
$BV$ functions. However it is the most natural way to introduce our framework. Recall a definition of $BV$ functions:
\begin{defi}
A function $f \in L^1(M)$ has \textbf{bounded variation} in $M$ if 
$$ \sup\{ \int_M f \; \div  \phi \;dx \; | \; \phi \in C_c^1(M,\ms{R}^n), |\phi|_{\infty} \leq 1\}< \infty.$$
In this case, $Df$ is defined by $\int_M f \; \div  \phi \;dx = -\int_M Df \; \phi \;dx$.
\end{defi}

\begin{defi}
We define $\text{Im}(M) \subset BV(M)$ such that for each function $f \in \text{Im}(M)$,
there exists a partition of $M$ in \lip domains $(U_i)_{i \in [0,k]}$ for an integer $k \geq 0$,
and the restriction $f_{|U_i}$ is \lip .
\end{defi}

\begin{rem}
The extension theorem of \lip function in $\ms{R}^n$ enables us to consider that on each $U_i$, $f_{| U_i}$ 
is the restriction of a \lip function defined on $\ms{R}^n$.  
\end{rem}

\noindent
On the definition of a \lip domain $U$: we use here (to shorten the previous definition) a large acceptation of 
\lip domains which can be found in chapter 2 of \cite{1002.49029}. Namely, $U$ is a lip domain if
there exists a \lip open set $\Omega$ such that $ \Omega \subset U \subset \bar \Omega$. 
In the proof of the derivation lemma \ref{pol},
we give the classical definition of \lip open set that we use above.
In a nutshell, an open set is \lip if for every point of the boundary there exists an affine basis of $\ms{R}^n$
in which we can describe the boundary of the open set as the graph of a \lip function on $\ms{R}^{n-1}$.
We chose to deal with \lip domains because it makes sense in the context of application to images.

\begin{example}
The most simple example is a piecewise constant function, $f = \sum_{i=1}^k a_i \mt{1}_{U_i}$ with $a_i \in \ms{R}$.
\end{example}

\begin{rem}
Our framework does not allow us to treat the discontinuities along a cusp, but we can deal with
the corners respecting the \lip condition. 
\end{rem}

\vspace{0.3cm}
\noindent
Let $f \in \text{Im}(M)$, we denote by $J_f$ the set of the jump part of $f$.
As a $BV(M)$ function, we can write the distributional derivative of $f$: 
$Df = \nabla f + D^c f + j(f)(x) \mathcal{H}^{n-1} \llcorner J_f$.
$\nabla f$ is the absolutely continuous part of the distributional derivative with respect to the Lebesgue measure and
$D_c$ is the Cantor part of the derivative. In other words, with the classical notations $j(f)(x)=(f^+(x)-f^-(x)) \nu_f(x)$, 
where $(f^+,f^-,\nu_f): J_f \mapsto \ms{R}^n \times \ms{R}^n \times \mt{S}^{n-1}$ is a Borel function. The functions $f^+$
and $f^-$ are respectively defined as $f^+(x) = \lim_{t \mapsto 0^+} f(x+t\nu_f(x))$ and $f^-(x) = \lim_{t \mapsto 0^-} f(x+t\nu_f(x))$.
Naturally, $j(f)$ does not depend on the choice of the representation of $\nu_f$, in fact $j(f)$ is homogeneous to the gradient. 
See for reference \cite{0909.49001} or \cite{MR1857292}.
In our case, the Cantor part is null from the definition. 
\\
We then write for $f \in \text{Im}(M)$, 

\begin{equation} \label{derivee}
Df = \nabla f + j(f)(x) \mathcal{H}^{n-1} \llcorner J_f.
\end{equation}

\subsection{The space of deformations}

\vspace{0.3cm}
\noindent
We denote by $V_M,<,>_V$ a Hilbert space of square integrable vector fields on $M$,
which can be continuously injected in $(\chi_0^p(M),\| . \|_{p,\infty})$, 
the vector space of $C^p$ with $p \geq 1$ vector fields which vanish on $\partial M$.
Hence, there exists a constant $c_V$ such that for all $v \in V$:
$$ \| v \|_{p,\infty} \leq c_V \|v \|_V.$$
Hence this Hilbert space is also a RKHS (Reproducing Kernel Hilbert Space), and we denote by
$k_V(x,.)\alpha$ the unique element of $H$ which verifies for all $v \in H$:
$<v(x),\alpha>=<k_V(x,.)\alpha,v>_H$, where $<,>$ is the euclidean scalar product
and $\alpha$ a vector in $\ms{R}^n$. This will enable an action on the support $M$.

\vspace{0.3cm}
\noindent
We denote by $S,<,>_S$ a Hilbert space of square integrable vector fields on $\ms{R}$,
as above. We denote by $k_S(x,.)$ its reproducing kernel.
This will enable the action on the level set of the functions.

\vspace{0.3cm}

\noindent
Through the following paragraph, we recall the well-known properties on the 
flow of such vector fields and its control. Most of them can be found in chapter $2$ of \cite{JoanPhD},
and are elementary applications of Gronwall inequalities. (See Appendix B in \cite{0902.35002})

\vspace{0.3cm}

\noindent
Let $v \in L^2([0,1],V)$, then with \cite{0855.57035} the flow is defined:
\begin{eqnarray} \label{flow}
\partial_t \phi_{0,t}^v &=& v_t \circ \phi_{0,t}^v, \\
\phi_0 &=& Id.
\end{eqnarray}
For all time $t \in [0,1]$, $\phi_{0,t}^v$ is a $C^1$ diffeomorphism of $M$ and the application 
$t \mapsto d_x \phi_{0,t}^v$ is continuous and solution of the equation:
\begin{equation} \label{deriv_diffeo}
d_x \phi_{0,t}^v = Id + \int_0^t d_{\phi_{0,s}^v(x)}v_s.d_x \phi_s^v ds.
\end{equation}

\vspace{0.3cm}

\noindent
We dispose of the following controls, with respect to the vector fields; let $u$ and $v$ be two vector fields
in $\in L^2([0,1],V)$ and $T \leq 1$:
\begin{eqnarray} \label{controle_champ}
\| \phi^u_{0,t}-\phi^v_{0,t} \|_ {\infty} &\leq & c_V \| v - u\|_{L^1[0,T]} \exp(c_V \| v \|_{L^1[0,T]}), \\
\| d\phi^u_{0,t} - d\phi^v_{0,t} \|_ {\infty} &\leq & 2 c_V \| v - u\|_{L^1[0,T]} \exp(c_V \| v \|_{L^1[0,T]}).
\end{eqnarray}

\vspace{0.3cm}
\noindent
And we have controls with respect to the time, with $[s,t] \subset [0,T]$:
\begin{eqnarray} \label{controle_temps}
\| \phi_{0,t}^v - \phi_{0,s}^v \|_{\infty} &\leq & \int_s^t \| v_r \|_{\infty} dr \leq c_V \int_s^t \| v_r \|_{V} dr,\\
\| \phi_{0,t}^v - \phi_{0,s}^v \|_{\infty} &\leq & c_V \sqrt{| s-t |} \| v \|_{L^2},\\
\| d\phi_{0,t}^v - d\phi_{0,s}^v \|_{\infty} &\leq & C \exp (C' \sqrt{T} \|v\|_{L^2[0,T]}) \int_s^t \| v_r \|_{V} dr.
\end{eqnarray}
\noindent
with the constants $C$ and $C'$ depending only on $c_V$.
Obviously these results are valid if $S$ replaces $V$.
In this case, we write $\eta_{0,t}$ for the flow generated by $s_t$. With the group relation for the flow, 
$\eta_{t,u} \circ \eta_{s,t} = \eta_{s,u}$.

\noindent
The group we consider is the product group of all the diffeomorphisms
we can obtain through the flow of $u \in L^2([0,1],V \times S)$. 

\vspace{0.5cm}

We aim to minimize the following quantity, with $\mu$ the Lebesgue measure:
\begin{equation} \label{J}
\mathcal{J}(v_t,s_t) = \frac{\lambda}{2} \int_{0}^1\|v_t\|_{V}^2dt + \frac{\beta}{2}\int_{0}^1\|s_t\|_{S}^2dt + \int_M |\eta_{0,1} \circ I_0 \circ \phi_{0,1}^{-1}(u)-I_{targ}(u)|^2d\mu(u) ,
\end{equation}
with $\eta_{0,t}$ the flow associated to $s_t$.
Remark that the metric we place on the product groups $V \times S$ is the product of the metric on each group
which is represented by the first two terms in \eqref{J}.
The functions $I_0,I_{targ}$ lie in $\text{Im}(M)$. In one section below, we prove classically that there exists 
at least one solution and we derive
the geodesic equations which give the form of the initial momentum.

\section{Derivation lemma} \label{lemmechap}
This derivation lemma may be useful in many situations where discontinuities arise.
Consider for exemple two \lip open sets $U$ and $V$. One may want to deform one of these open sets while
the second remains unchanged (figure below). The basic case is the following: 
$$J_t=\int_V \chi_U \circ \phi_t^{-1} dx = \mu(V \cap \phi_t(U)),$$
with $\mu$ the Lebesgue measure. We answer to the differentiation of $J_t$, we obtain a sort of Stokes formula 
with a perturbation term. We discuss below a more general formula to apply in our context. 
The final result is the proposition \ref{lip1}:

\begin{lm} \label{lip1}
Let $U,V$ two bounded \lip domains of $\mt{R}^n$. 
Let $X$ a \lip vector field on $\ms{R}^n$ and $\phi_t$ the associated flow. 
Finally, let $g$ and $f$ \lip real functions on $\mt{R}^n$.
Consider the following quantity depending on $t$,
$$ J_t = \int_{\phi_t(U)} f \circ \phi_{t}^{-1} g \mt{1}_V d\mu,$$
where $d\mu$ is the Lebesgue measure, then
\begin{equation}
\partial_{t|t=0^+} J_t = \int_{U} -<\nabla f,X> g \mt{1}_V d\mu + \int_{\partial U} <X,n> fg  \tilde{\mt{1}}_V(X) d\mu_{|\partial U}.
\end{equation}
with $ \tilde{\mt{1}}_V(X)(y)= \lim_{\epsilon \mapsto 0^+} \mt{1}_{\bar V}(y+\epsilon X) $, if the limit exists, $0$ elsewhere.
And we denote by $d\mu_{|\partial U}$ the measure on $ \partial U$ and $n$ the outer unit normal.
\end{lm}

\vspace{0.2cm}

\begin{center}
\begin{figure} 
\includegraphics{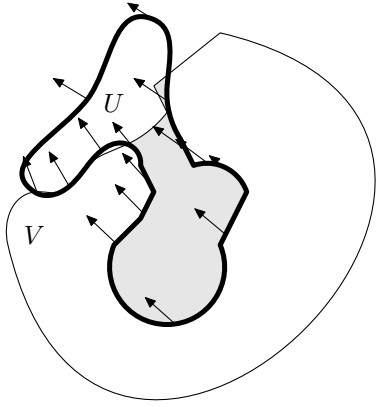}    
\caption{Evolution of the area between two Lipschitz curves. (The arrows represent $X$ along
the boundary of $U$)}
\end{figure}
\end{center}

\noindent
As a corollary, we deduce:

\begin{corr}
We have,
$$\partial_{t|t=0^+} \mu(V \cap \phi_t(U)) = \int_{\partial U} <X,n> \tilde{\mt{1}}_V(X) d\mu_{|\partial U},$$
with $ \tilde{\mt{1}}_V(X)(y)= \lim_{\epsilon \mapsto 0^+} \mt{1}_{\bar V}(y+\epsilon X) $, if the limit exists, $0$ elsewhere.
And we denote by $d\mu_{|\partial U}$ the measure on $ \partial U$ and $n$ the outer unit normal.
\end{corr}

\noindent 
In this case, the derivation formula is a Stokes' formula in which one takes only into account the deformation viewed in $V$.

\vspace{0.2cm}

\noindent
Below is a figure to illustrate the lemma:
\begin{rem}
We could generalize the lemma to finite intersection of \lip domains, with the same scheme of the proof
developed above. We gain hence generality which seems to be very natural for concrete applications.
\end{rem}

\noindent
This generalization for \lip domains is sufficient for the application we aim, and this application is
presented in the paragraph below to derive the geodesic equations. Hence we present the corollary we use in the next 
paragraph.

\begin{theo}
Let $(f,g) \in \text{Im}(M)^2$, $X$ a \lip vector field on $\ms{R}^n$ and $\phi_t$ the associated flow.
$$ J_t = \int_M f \circ \phi_t^{-1} (x) g(x) d\mu(x),$$
then the derivation of $J_t$ is:
\begin{equation} \label{lemma1}
\partial_{t|t=0+} J_t = \int_{M} -<\nabla f,X> g dx - \int (f_+-f_-)  \tilde{g} <\nu_f,X> d\mathcal{H}^{n-1}.
\end{equation}
with $\tilde{g}_X(x):=\lim_{t \mapsto 0^+} g(\phi_t(x))$ if the limit exists and if not, $\tilde g_X(x)=0$.
\end{theo}

\vspace{0.3cm}

\noindent
\emph{Proof:}
\noindent
Writing $f$ as $f=\sum_{i=1}^{k} f \mt{1}_{x \in U_i}$ where $(U_i)_{i=1,\ldots,n}$ is the partition in domains associated to $f$,
and using the same expression for $g$, by linearity of integration, we fall in the case of the proposition \ref{lemmechap}.
$\square$

\vspace{0.3cm}

\noindent
A last remark on the formulation of the lemma, we can rewrite the equation \eqref{lemma1} in a more compact form:
$$\partial_{t|t=0+} J_t = - \int_{M} \langle Df,X \rangle \tilde g, $$ 
with $Df$ the notation for the derivative for $SBV$ function and $\tilde g$ is the function defined above. Remark 
that $\mu \; a.e.$ $\tilde g= g$, these two functions differ on $J_f$.

\section{Minimizing the energy functional}
The existence of geodesics is a classical fact, but in this framework the derivation of the geodesic equations did
not appear to the author in the existing literature. 
With the metric introduced above, a geodesic in the product space is a product of geodesics.
We chose to understand the two geodesics separately for technical reasons.
We could also have described the geodesics in the product space $V \times S$,
this point of view will be detailed in the Hamiltonian formulation of the equations.

\subsection{Existence and equations of geodesics}
\begin{theo}
Let $(I_0,I_{tar}) \in \text{Im}(M)^2$, we consider the functional $\mathcal{J}$ on $H=L^2([0,1],V \times S)$ defined in \eqref{J}.
There exists $(v,s) \in H$ such that $\mathcal{J}(v,s)= \min_{(v,s) \in H} \mathcal{J}(v,s)$.
For such a minimizer, there exists $(p_a,p_b,p_c) \in L^2(M,\ms{R}^n) \times L^2(J_{I_0},\ms{R}^n) \times L^2(M,\ms{R})$ such that:
\begin{eqnarray}
\beta s_t = \int_M  p_c(y) d[\eta_{t,1}]_{I_t^s(y)}  k_S(I^s_t(y),.) d\mu(y), \label{contraste} \\
\lambda v_t = \int_M k_V( \phi_{0,t} (x) ,.) [d \phi_{0,t}]_x^{-1*} (p_a(x)) d\mu(x)  \nonumber \\
+ \int_{J_{I_0}} k_V( \phi_{0,t} (x) ,.) [d \phi_{0,t}]_x^{-1*} (p_b(x)) d\mu_{ | J_{I_0}} (x), \label{generale}
\end{eqnarray}
with:
\begin{eqnarray*}
I_t^s=\eta_{0,t} \circ I_0 \circ \phi_{1,0}, \\
I_t^v=\eta_{0,1} \circ I_0 \circ \phi_{t,0},
\end{eqnarray*}
and $J_{I_0}$ the jump set of $I_0$.
More precisely for the $(p_a,p_c)$ we show, we have the equation:
\begin{equation} \label{relation}
p_a(x) = -\Jac (\phi_{0,1}(x)) \nabla_{|x}I_0^v \; p_c(x).
\end{equation}
\end{theo}

\vspace{0.3cm}

\noindent
\emph{Proof:}

\vspace{0.2cm}

\noindent
On the space $H$, the strong closed balls are compact for the weak topology.
The functional $\mathcal{J}$ is lower semi-continuous, so we obtain the existence of a minimizer.
As reference for the weak topology \cite{00865216}.  
We find here \cite{trouve05:_local_geomet_defor_templ} a proof that the flow is continuous for the weak topology, 
the main point to prove the semi-continuity: if $(u_n,s_n) \rightharpoonup (u,v)$ in $H$, 
then $(\phi_{0,1},\eta_{0,1}) \mapsto (\phi_{0,1},\eta_{0,1})$.

\vspace{0.2cm}
\noindent
We first differentiate w.r.t. the vector field $s \in L^2([0,T],S)$, we denote by $\tilde s$ a perturbation
of $s$. Using the lemma \ref{depder} in appendix, we write:
$\partial_{\tilde s} \eta_{0,1}(x) = \int_0^1 [d\eta_{t,1}]_{|\eta_{0,t}(x)}\tilde s_t(\eta_{0,t}(x)) dt$.  
We already introduce the kernel:
$$ \int_0^{1}[ \beta <s_t,\tilde s_t>_S + \int_M 2(I_1-I_{targ}) 
[d\eta_{t,1}]_{|\eta_{0,t}(I_0^s(y))} <k_S(I_t^s(y),.),\tilde s_t>_S d\mu(y) ]dt = 0,$$ 
it leads to:
$$ \beta s_t + \int_M 2(I_1(y)-I_{targ}(y)) [d\eta_{t,1}]_{|I_t^s(y)} k_S(I_t^s(y),.) d\mu(y) = 0.$$ 
With the notation $p_c = -2(I_1-I_{targ})$, we have the first equation announced.

\vspace{0.2cm}
\noindent
For the second equation, we need the derivation lemma detailed in section \ref{lemmechap}. In order to use the lemma,
we first need to develop the attachment term:

$$ \int_M |\eta_{0,1} \circ I_0 \circ \phi_{0,1}^{-1}(u)-I_{targ}(u)|^2d\mu(u)= \int_M (I_1^v)^2 - 2I_1^vI_{targ} + I_{targ}^2 d\mu,$$
Now, only the first two terms are involved in the derivation, and we apply the lemma to these two terms. 
(Actually the lemma is necessary only for the second term.) 
\\
Again, we have with the lemma \ref{depder}:
$$ V_1:= \partial_{\epsilon} \phi_{1,0}(x) = - \int_0^1 d (\phi_{t,0})_{\phi_{1,t} (x)} (\tilde v [ \phi_{1,t} (x) ]) dt.$$

\noindent
We consider the semi-derivation of (\ref{J}) at the minimum with respect to the displacement field $v$, we use the notations
of $SBV$ functions for the derivatives:

\begin{eqnarray*} \label{pseudoder}
&\phantom{1}& \lambda \int_0^1 <v_t,\tilde v_t> dt + \int D[(I^v_1)^2] V_1 -2 \int DI^v_1 V_1 \tilde I_{targ}^{\tilde v} =0, \\
&\phantom{1}& \lambda \int_0^1 <v_t,\tilde v_t> dt + \int_M 2 <\nabla I_1^v,V_1> (I_1^v - I_{targ}) d\mu + \\
&\phantom{1}& \phantom{dddddddddddd}\int_{J_{I_1^v}} \left( j([I_1^v]^2)-2j(I_1^v) \right) V_1 \tilde I_{targ} d\mathcal{H}^{n-1}=0.
\end{eqnarray*}
As $(I^v_1)^2$ and $I^v_1$ have the same discontinuity set, the second integration is only over $J_{I_1^v}$. 

\noindent
We apply a version of the central lemma in \cite{hamcurves} which is detailed in appendix (see lemma \ref{lmc}).

\begin{eqnarray*}
g: L^2([0,1],V) &\mapsto & L^2([0,1],V) \times L^2(M,\ms{R}^n) \times L^2(J_{I^v_1},\ms{R}^n) \\
	\tilde v &\mapsto & (\tilde v,V_1,{V_1}_{|J_{I_v^1}}), \\            
\end{eqnarray*}

\noindent
We ensure that $B:=\{(v,2 \nabla I_1^v (I_1^v - I_{targ}), \left( j([I_1^v]^2)-2j(I_1^v) \right) \tilde I_{targ}), \tilde v \in L^2([0,1],V) \}$
is bounded. For each $\tilde v$, $|\tilde I_{targ}^{\tilde v}|_{\infty} \leq | I_{targ}|_{\infty}$.
(This assumption could be weakened.) Whence we get with the lemma \ref{lmc} the existence of $\tilde I_{targ} \in \overline{\Conv(B)}$
(we observe that the Lebesgue part of $\tilde I_{targ}$ is equal to $I_{targ}$, the modification is on the set $J_{I_1^v}$),
 such that:

\begin{eqnarray*}
&\phantom{0}&\lambda  \int_0^1 <v_t,\tilde v_t> dt + \; \int D[(I^v_1)^2] V_1 -2 \int DI^v_1 V_1 \tilde I_{targ} = 0, \\
&\phantom{0}&\lambda  \int_0^1 <v_t,\tilde v_t> dt + \; \int_M 2 \nabla I_1^v (I_1^v - I_{targ}) d\mu + 
\int_{J_{I_1^v}} j([I_1^v]^2)V_1-\int_{J_{I_1^v}} 2j(I_1^v)\tilde I_{targ}V_1=0. 
\end{eqnarray*}

\noindent
Now, we aim to obtain the explicit geodesic equations by introducing the kernel, we denote by 
$A(\tilde v):=\int D[(I^v_1)^2] V_1 -2 \int DI^v_1 V_1 \tilde I_{targ}$ the pseudo derivative of the attachment term and
we denote also:
$$ \tilde \Delta (x) :=  j((I_0)^2) \circ \phi_{1,0}(x) - 2 j(I_0) \circ \phi_{1,0}(x) \tilde I_{targ}(x),$$ 
which defines a normal vector field on $J_{I_1^v}$. 

\begin{eqnarray*}
A(\tilde v) = \int_{0}^{1} - \int_M 2 (I_1-I_{targ})(y)  < k_V( \phi_{1,t} (y) ,.) d(\phi_{t,0})_{\phi_{1,t} (y)}^{*} (\nabla_{| \phi_{1,0}(y)} I_0^v),\tilde v_t>_V d\mu(y)  \\
- \int_{\phi_{0,1} (J_{I_0})}   < k_V( \phi_{1,t} (y) ,.) d(\phi_{t,0})_{\phi_{1,t} (y)}^{*} (\tilde \Delta(y)),\tilde v_t(y)>_V d\mu_{| \phi_{0,1}(J_{I_0})} (y) dt.
\end{eqnarray*}

\noindent
With the change of variable $x=\phi_{1,0}(y)$,
\begin{eqnarray} \label{der}
A (\tilde v) = \int_{0}^{1} - \int_M  2 (I_1-I_{targ})(\phi_{0,1}(x)) \Jac (\phi_{0,1}(x)) < k_V( \phi_{0,t} (x) ,.) d (\phi_{t,0})_{\phi_{0,t} (x)}^{*} (\nabla_{| x}I_0^v),\tilde v_t>_V d\mu \nonumber \\ 
- \int_{J_{I_0}}   \frac{\Jac (\phi_{0,1}(x))}{|d\phi_{0,1} (n_x)|} < k_V( \phi_{0,t} (x) ,.) d(\phi_{t,0})_{\phi_{0,t} (x)}^{*} (\tilde \Delta(\phi_{0,1}(x)),\tilde v_t>_V d\mu_{| J_{I_0}} dt,
\end{eqnarray}
with $n_x$ a normal unit vector to $J_{I_0}$ in $x \in J_{I_0}$. Note that in the second term, the change of variable acts on the hypersurface 
$J_{I_0}$. This explains the term  $\frac{\Jac (\phi_{0,1}(x))}{|d\phi_{0,1} (n_x)|}$ which corresponds to the Jacobian term 
for the smooth part.

\vspace{0.3cm}
\noindent
We are done,
\begin{eqnarray*}
\lambda v_t = \int_M 2 \Delta(\phi_{0,1}(x)) \Jac (\phi_{0,1}(x)) k_V( \phi_{0,t} (x) ,.) 
d (\phi_{0,t}^{-1})_{\phi_{0,t} (x)}^{*} (\nabla_{| x} I_0^v) d\mu \\ 
+ \int_{J_{I_0}}  \frac{\Jac (\phi_{0,1}(x))}{|d\phi_{0,1} (n_x)|}
k_V( \phi_{0,t} (x) ,.) d(\phi_{0,t}^{-1})_{\phi_{0,t} (x)}^{*} ( \tilde \Delta(\phi_{0,1}(x)) ) d\mu_{ |J_{I_0}}. 
\end{eqnarray*}

\noindent
With $p_a(x)= 2\Delta(\phi_{0,1}(x)) \Jac (\phi_{0,1}(x))\nabla_{| x}(I_0^v)$ and 
$p_b(x) = \tilde \Delta(\phi_{0,1}(x)) \frac{\Jac (\phi_{0,1}(x))}{|d\phi_{0,1} (n_x)|}$,	
we have the geodesic equations.
$\square$

\vspace{0.3cm}
\noindent
These geodesic equations are a necessary condition for optimality. In the next paragraph, we show that if $(p_a,p_b,p_c)$
is given, we can reconstruct the geodesics.

\subsection{Reconstruction of geodesics with the initial momentum}

We first demonstrate that if a vector field is a solution to the geodesic equations, then the norm is constant in time.
\begin{prop} \label{vitesseconstante}
\textbf{Constant speed curves in vector fields spaces}
\\
If a vector field $s_t$ is a solution of equation \eqref{contraste} 
and the kernel is differentiable then $\| s_t \|^2$ is constant.
\\
If a vector field $v_t$ is a solution of equation \eqref{generale} 
and the kernel is differentiable then $\| v_t \|^2$ is constant.
\end{prop}

\vspace{0.3cm}

\noindent
\emph{Proof:}
\\
\noindent
We prove the first point:
$$\| s_t \|_S^2 = \int_M \int_{M} p(y') d[\eta_{t,1}]_{I_t^s(y')} k_S(I^s_t(y'),I^s_t(y)) p(y) d[\eta_{t,1}]_{I_t^s(y)} d\mu(y')d\mu(y).$$

\noindent
Remark that a.e. $\partial_t (d[\eta_{t,1}]_{I_t^s(y)}) = -d[s_t]_{I_t^s(y)} d[\eta_{t,1}]_{I_t^s(y)}$.
This equation is obtained by a derivation of the group relation:
$ d[\eta_{0,1}]_{I_0^s} = d[\eta_{t,1} \circ \eta_{0,t}]_{I_0^{s}}$, and with the derivation of the equation \eqref{contraste}:
$$d[s_t]_{x}=\int_M d[\eta_{t,1}]^*_{I_t^s(y)}p_c(y) \partial_1 k_S(x,I^s_t(y))dy.$$

\noindent
As $ds_t \in L^1([0,1])$ then the equation \eqref{deriv_diffeo} proves that $d[\eta_{t,1}]_{I_t^s(.)}$
is absolutely continuous. As the space of absolutely continuous functions is an algebra,
$\| s_t \|_S^2$  is also absolutely continuous.
To obtain the result, it suffices to prove that the derivate vanishes a.e. 

\begin{eqnarray*}
\partial_t \| s_t \|^2 = -\int_M \int_{M} p(y') ds_t(I_t^s(y')) d[\eta_{t,1}]_{I_t^s(y')} k_S(I^s_t(y')),I^s_t(y))) p(y) d[\eta_{t,1}]_{I_t^s(y)} d\mu(y')d\mu(y)\\
+\int_M \int_{M} p(y') d[\eta_{t,1}]_{I_t^s(y')} \partial_1 k_S(I^s_t(y')),I^s_t(y)))s_t(I^s_t(y')) p(y) d[\eta_{t,1}]_{I_t^s(y)} d\mu(y')d\mu(y) = 0.
\end{eqnarray*}

\vspace{0.3cm}
\noindent
The second point is very similar. 
We underline that the equation \eqref{generale} is a particular case of the following, with a measure $\nu$ which has a Lebesgue part 
and a singular part on the set $J_{I_0}$ of discontinuities of the function $I_0$. We also define:
$$p_t(x) =  (d[\phi_{0,t}]^*_x)^{-1}(p_a(x) \mt{1}_{x \notin J_{I_0}} + p_b(x) \mt{1}_{x \in J_{I_0}}),$$
By the definition,
$$\| v_t \|^2 = \int \int p_t(x) k_V( \phi_{0,t} (x) , \phi_{0,t} (y) ) p_t(y) d\nu(x) d\nu(y).$$
Remark that
$ \partial_t p_t(x) = - dv_t^* \circ \phi_{0,t}(x) p_t(x)$,
and we differentiate:
\begin{eqnarray*}
\partial_t \| v_t \|^2 = -\int \int \int p_t(x) \partial k_V(\phi_{0,t}(x),\phi_{0,t}(z)) p_t(z) k_V( \phi_{0,t} (x) , \phi_{0,t} (y) ) p_t(y) d\nu(x) d\nu(y) d\nu(z) \\
+ \int \int \int p_t(x) \partial k_V(\phi_{0,t}(x),\phi_{0,t}(y)) p_t(y)  k_V( \phi_{0,t} (x) , \phi_{0,t} (z)) p_t(z) d\nu(x) d\nu(y) d\nu(z) =0.
\end{eqnarray*}
$\square$

\vspace{0.3cm}

\noindent
This proposition is crucial to establish that the geodesics are defined for all time.
Namely, we answer to existence and uniqueness of solutions to (the set $J_{I_0}$ 
but could be much more general than the discontinuity set of a function in $\text{Im}(M)$):

\begin{eqnarray} \label{sys}
\eta_{0,t} &=& Id + \int_0^t s_u \circ \eta_u du, \nonumber \\
\beta s_t(.) &=& \int_M  p_c(y) d[\eta_{t,0}]_{I_t^s(y)}  k_S(I_t^s(y)),.) d\mu(y), \nonumber \\
\phi_{0,t} &=& Id + \int_0^t v_u \circ \phi_u du, \nonumber \\
\lambda v_t(.) &=& \int_M k_V(., \phi_{0,t} (x)) [d \phi_{0,t}]_x^{-1*} (p_a(x)) d\mu(x) \nonumber\\
&+& \int_{J_{I_0}} k_V(., \phi_{0,t} (x)) [d \phi_{0,t}]_x^{-1*} (p_b(x)) d\mu_{ | J_{I_0}} (x) dt.
\end{eqnarray}

\vspace{0.3cm}
\noindent
On purpose, this system of equations is decoupled in $v$ and $s$. The proof of the next proposition treats both cases
in the same time but it could be separated.

\begin{prop} \label{temps}
For $T$ sufficiently small, the system of equations \eqref{sys} with 
$$(p_a,p_b,p_c) \in L^1(M,\ms{R}^n) \times L^1(J_{I_0},\ms{R}^n) \times L^1(M,\ms{R})$$ has a unique solution 
if the kernel is differentiable and its first derivative is Lipschitz.
\end{prop}

\vspace{0.3cm}

\noindent
\emph{Proof:}
\\
\noindent
We aim to apply the fixed point theorem on the Banach space $L^2([0,T],V \times S)$.
We estimate the \lip coefficient of the following application:

\begin{eqnarray} \label{Xi}
\Xi : L^2([0,T],V \times S)  & \mapsto & L^2([0,T],V \times S)  \nonumber \\
(v,s) & \mapsto & (\xi(v),\xi(s)),
\end{eqnarray}
with
\begin{eqnarray*}
\xi(s)_t = \int_M  p_c(x) d[\eta_{t,0}]_{\tilde I_t(x)}  k_S(\tilde I_t(x)),.) d\mu(x), & \nonumber \\
\xi(v)_t = \int_M k_V(., \phi_{0,t} (x)) [d \phi_{0,t}]_x^{-1*} (p_a(x)) d\mu(x) & \nonumber\\
+ \int_{J_{I_0}} k_V(., \phi_{0,t} (x)) [d \phi_{0,t}]_x^{-1*} (p_b(x)) d\mu_{ | J_{I_0}} (x) dt.& 
\end{eqnarray*}
\noindent
For the space $L^2([0,T],V)$, if we have:
$$\|\xi(v)_t-\xi(u)_t\|^2 \leq M \| v - u\|_{L^1[0,T]},$$
the result is then proven with Cauchy-Schwarz inequality:
$$\|\xi(v)-\xi(u)\|_{L^2[0,T]} \leq \sqrt{MT}\| v - u\|_{L^2[0,T]}.$$
This can be obtained with:
\begin{eqnarray*}
\|\xi(v)_t-\xi(u)_t\|^2 &=& <\xi(v)_t,\xi(v)_t-\xi(u)_t> - <\xi(u)_t,\xi(v)_t-\xi(u)_t> , \\
\|\xi(v)_t-\xi(u)_t\|^2 &\leq & 2 \max (|<\xi(v)_t,\xi(v)_t-\xi(u)_t>|,|<\xi(u)_t,\xi(v)_t-\xi(u)_t>|).
\end{eqnarray*}
\noindent
For one of the two terms in the equation above:
\begin{eqnarray*}<\xi(v)_t,\xi(v)_t-\xi(u)_t> &=& \int \int [d \phi_{0,t}^v]_x^{-1*} (p_a(x)) [k(\phi_{0,t}^v (x),\phi_{0,t}^v (y)) [d \phi_{0,t}^v]_y^{-1*} (p_a(y)) \\
&-& k(\phi_{0,t}^v (x),\phi_{0,t}^u(y)) [d \phi_{0,t}^u]_y^{-1*} (p_a(y))] d\mu(y) d\mu(x).
\end{eqnarray*}
On the unit ball of $L^2([0,T],V \times S)$ denoted by $B$, and with the inequality \eqref{controle_champ}, we control the diffeomorphisms:
\begin{eqnarray} \label{controle_champ2}
\| \phi^u_{0,t}-\phi^v_{0,t} \|_ {\infty} &\leq & c_V \| v - u\|_{L^1[0,T]} \exp(c_V), \\ \nonumber
\| d\phi^u_{0,t} - d\phi^v_{0,t} \|_ {\infty} &\leq & 2 c_V \| v - u\|_{L^1[0,T]} \exp(c_V). \nonumber
\end{eqnarray}
With the triangle inequality, we get:
\begin{eqnarray*} 
|<\xi(v)_t,\xi(v)_t-\xi(u)_t>| & \leq & \int \int  [d \phi_{0,t}^v]_x^{-1*} (p_a(x)) |[k(\phi_{0,t}^v (x),\phi_{0,t}^v (y)) [d \phi_{0,t}^v]_y^{-1*} (p_a(y)) \\
&-& k(\phi_{0,t}^v (x),\phi_{0,t}^v(y)) [d \phi_{0,t}^u]_y^{-1*} (p_a(y))]| \\
&+& |[k(\phi_{0,t}^v (x),\phi_{0,t}^v (y)) [d \phi_{0,t}^u]_y^{-1*} (p_a(y)) \\
&-& k(\phi_{0,t}^v(x),\phi_{0,t}^u(y)) [d \phi_{0,t}^u]_y^{-1*} (p_a(y))]| d\mu(x) d\mu(y).
\end{eqnarray*}

\noindent
On the unit ball $B$, we have:
\begin{eqnarray*}
&\| d\phi^u_{0,t} \|_{\infty} &\leq   1+2c_V \exp(c_V), \\
&\| \phi^u_{0,t}-Id \|_{\infty} &\leq  c_V.
\end{eqnarray*}
\noindent
Let $M_k \in \ms{R}$ a bound for the kernel and its first derivative on the unit ball $B$.
Such a constant exists thanks to the hypothesis on the kernel and its first derivative.

\noindent
A bound for the first term can be found with the second inequality of \eqref{controle_champ2}:
$$ 2 c_V \| v - u\|_{L^1[0,T]} \exp(c_V)M_k (1+2c_V \exp(c_V)) \| p_a\| \|p_b\|,$$
the second term is controlled with the first inequality of \eqref{controle_champ2} 
with the \lip hypothesis on the kernel:
$$  c_V \| v - u\|_{L^1[0,T]} \exp(c_V)M_k (1+2c_V \exp(c_V)) \| p_a\| \|p_b\|.$$
\noindent
Finally we get,
$$\| \xi(v)_t-\xi(u)_t\|^2 \leq  6 \; c_V \| v - u\|_{L^1[0,T]} \exp(c_V)M_k (1+2c_V \exp(c_V)) \| p_a\| \|p_b\|.$$
We have now concluded for the first component of the application $\Xi$.
For the second term, the proof is essentially the same, we do not give the details.
$\square$

\noindent
We have proven that there exists $T>0$ such that we have existence and uniqueness to the system \eqref{sys},
we prove now that the solutions are non-exploding i.e. we can choose $T=+\infty$ in the last proposition.
This property shows that the associated riemannian manifold of infinite dimension is complete, since
the exponential map is defined for all time. Without the hypothesis on the kernel, we can find simple counter-examples to 
this fact.

\begin{prop} 
The solution proposition \ref{temps} is defined for all time.
\end{prop}

\noindent 
\emph{Proof:}
\\
\noindent
Thanks to propostion \ref{vitesseconstante}, we know that the norm of the solution $u_t$ is constant in time, 
which will enable the extension for all time. 
Consider a maximal solution with interval of definition $[0,T]$ with $T < +\infty$,
then with the inequalities from \eqref{controle_temps} and after, we define the limit $\lim_{t \mapsto T} \phi_{0,t}:=\phi_{0,T}$,
since for all $x$, $\phi_{0,t}(x)$ is a Cauchy sequence. This is the same for $\lim_{t \mapsto T} d\phi_{0,t}(x)$.
This limit is also a diffeomorphism, since we can define the limit of the inverse as well.
The proof is the same to extend $\eta_t$ for all time.
\\
We can then apply the proposition (existence for small time) to the current image $I_t$ instead of $I_0$,
we obtain diffeomorphisms $\tilde \phi_{0,s}$ and $\tilde \eta_{0,s}$ in a neighborhood of $0$, $[0,\epsilon]$.
Composing with $\phi_{0,T}$ and $\eta_{0,T}$, we extend the maximal solution on $[0,T+\epsilon]$.
This is a contradiction.
$\square$

\vspace{0.3cm}
\noindent
We decoupled the equations in $s$ and $v$ to give a simple proof of the existence in all time of the flow.
The formulation of \eqref{sys} implies the following formulation, which is the first step to understand the 
weak Hamiltonian formulation. If we have the system \eqref{sys} and the relation \eqref{relation},
through the change of variable $u=\phi_{0,t}(x)$, we get easily:

\begin{eqnarray} \label{sys2}
\eta_{0,t} &=& Id + \int_0^t s_u \circ \eta_u du, \nonumber \\
\beta s_t(.) &=& -\int_M  P_t(x) k_S(I_t(x)),.) d\mu(x), \nonumber \\
\phi_{0,t} &=& Id + \int_0^t v_u \circ \phi_u du, \nonumber \\
\lambda v_t(.) &=& \int_M k_V(., u) P_t(u) \nabla_{|u} I_t d\mu(x) , \nonumber\\
&+& \int_{J_{I_0}} k_V(., \phi_{0,t} (x)) [d \phi_{0,t}]_x^{-1*} (p_b(x)) d\mu_{ | J_{I_0}} (x) dt,
\end{eqnarray}
with,
\begin{eqnarray*}
P_t(x) &=& \Jac (\phi_{t,0}) d[\eta_{t,0}]_{I_t^s(x)} P_0 \circ \phi_{t,0}, \\
P_0(x) &=& -p_c(\phi_{0,1}(x))d[\eta_{0,1}]_{I_0(x)} \Jac (\phi_{0,1}(x)).
\end{eqnarray*}

\section{Hamiltonian generalization}

In numerous papers on large deformation diffeomorphisms, the Hamiltonian framework arises. The simplest example
is probably the Landmark matching problem for which the geodesic equations and the Hamiltonian version of the evolution
are well known (\cite{vmty04},\cite{bb37524}).
Our goal is to provide an Hamiltonian interpretation of the initial variational problem. The main difference is
that we want to write Hamiltonian equations in an infinite dimensional space, which is roughly the space of images.
The first step was done in \cite{hamcurves} where Hamiltonian equations were written on the representations of 
closed curves. Our work generalizes this approach to the space of images.
We use the point of view of the optimal control theory (as it is developed in \cite{bb37524}) to formally introduce 
the Hamiltonian. Then, we state a weak Hamiltonian formulation of the equations obtained by the variational approach.
We prove uniqueness for the solutions to these equations. At the end of this section, we discuss the existence of the solutions
with the help of the existence of the solutions for the variational problem.
\\
In the whole section, we maintain our previous assumptions on the kernel
for the existence of solutions for all time.

\subsection{Weak formulation}

\vspace{0.3cm}
\noindent
In this paragraph, we slightly modify the approach in order to develop the idea of decomposing an image
in "more simple parts".
Let introduce the position variables. We consider in the following that the discontinuity boundary is 
a position variable. Instead of considering the function $I_t$ as the second position variable, 
we introduce a product space which can be projected on the space $\text{Im}(M)$.
Let $(U_1,\ldots,U_n)$ be a partition in \lip domain of $M$. We denote by $\Sigma_0 = \cup_{i=1}^{n} \partial U_i$ 
the union of the boundaries of the \lip domains. 
We consider the projection:
\begin{eqnarray}
p: \prod_{i=1}^{n} W^{1,\infty}(M) & \mapsto & \text{Im}(M) \\
	(I^i)_{i=1,\ldots ,n} & \mapsto & I=\sum_{i=1}^n I^i \mt{1}_{U_i}.
\end{eqnarray}	
Discontinuities give derivatives with a singular part. Maybe we could have treated 
this case adopting only the variable $I_t$, but we find the idea of decomposing an image 
into more simple parts rich enough to study the case.
Observe that $\Sigma_0$ is endowed with an important role in the definition of the projection: 
to write down a Hamiltonian system on the large space, we need to introduce the deformation of 
$\Sigma_0$ and the deformation of each function in the product space. We will derive the Hamiltonian equations
from this optimal control problem:
(the position variable is $Q$ and the control variable is $U$, $c(U)$ is the instantaneous cost function) 

\begin{eqnarray*}
Q &=& (Q_i)_{0 \leq i \leq  r} = (\Sigma,(I^i)_{1 \leq i \leq r}) \in L^1(\Sigma_0,M) \times W^{1,\infty}(M)^r, \\
U &=& (v,s) \in V \times S, \\
\dot Q &=& f(Q,U) = (v \circ Q^0,(-\langle \nabla Q^i,v \rangle + s(Q^i))_{1 \leq i \leq  r}), \\
c(U) &=& \frac{\lambda}{2}|v|_V^2 + \frac{\beta}{2}|s|_S^2.
\end{eqnarray*}

\noindent
The cotangent space $D$ of the position variable contains $F=L^{\infty}(\Sigma,\ms{R}^n) \times L^1(M,\ms{R})^r \subset D$.
We write the \textbf{formal} minimized Hamiltonian of the control system on the subspace $F$, with $P \in F$:

\begin{equation} \label{eq_ham}
H(P,Q) = \min_{U} \int_{\Sigma_0} \langle P_0(x), \dot Q_0(x) \rangle d\mu_{| \Sigma_0}(x) + \sum_{i=1}^r \int_M P^i \dot Q_i d\mu(x) - c(U).
\end{equation}

\noindent
Minimizing in $U$, we obtain optimality conditions in $(u,v)$ a minimizer such that for any perturbation $(\delta v,\delta s)$:

\begin{eqnarray*}
\lambda \langle v, \delta v \rangle&=& \int_{\Sigma_0} \langle P^0,\delta v \circ Q^0 \rangle d\mu_{| \Sigma_0}(x)  - \sum_{i=1}^r \int_M P^i \langle \nabla Q^i,\delta v \rangle d\mu(x), \\
\beta \langle s,\delta s \rangle &=& \sum_{i=1}^r \int_M P^i \delta s(I^i) d\mu(x).
\end{eqnarray*}

\noindent
Using the kernel, it can be rewritten,
\begin{eqnarray} \label{def_v_s}
\lambda v &=& \int_{\Sigma_0} k(Q^0(x),.)P^0 d\mu_{| \Sigma_0}(x) - \sum_{i=1}^r \int_M  k_V(x,.)P^i \nabla Q^i d\mu(x), \\
\beta s &=& \sum_{i=1}^r \int_M k_S(Q^i,.)P^i d\mu(x).
\end{eqnarray}

\noindent
We deduce the expression of the Hamiltonian,
\begin{eqnarray*}
H(P,Q) &=& \frac{1}{2\lambda}[\int_{\Sigma_0} \int_{\Sigma_0} P^0(x) k_V(Q^0(x),Q^0(y)) P^0(y) d\mu_{| \Sigma_0}(x)d\mu_{| \Sigma_0}(y) \\
&+& \int_M \int_M P^j(y) \nabla Q^j(y) k_S(y,x) \nabla Q^i(x) P^i(x) d\mu(x)d\mu(y) \\
&-& 2 \sum_{1 \leq i \leq r} \int_M \int_{\Sigma_0}  P^0(y)k_V(Q^0(y),x)P^i(x) \nabla Q^i(x) d\mu(x)d\mu_{| \Sigma_0}(y)] \\ 
&+& \frac{1}{2\beta} \sum_{1 \leq i,j \leq r} \int_M \int_M P^j(y) k_S(Q^j(y),Q^i(x)) P^i(x) d\mu(x)d\mu(y).
\end{eqnarray*}

\noindent
Now, we want to give a sense to the Hamiltonian equations, $\forall \; i \in [1,r]$:

\begin{eqnarray} \label{ham}
\dot Q^0_t  &=& \partial_{P^0} H(P_t,Q_t) (.), \\ \nonumber
\dot Q_t^i &=& \partial_{P^i} H(P_t,Q_t) \; \forall \; i \in [1,r], \\ \nonumber
\dot P_t^0 &=& -\partial_{Q^0} H(P_t,Q_t), \\ \nonumber
\dot P_t^i &=& - \partial_{Q^i}  H(P_t,Q_t) \; \forall \; i \in [1,r].\\ \nonumber
\end{eqnarray}

\noindent
These derivatives should be understood as distributions, 
for $\Psi \in C_0^{+\infty}(M,\ms{R})$ and $u \in C_0^{+\infty}(M,\ms{R}^n)$ and with the notation introduced in (\ref{def_v_s}),
$ \forall \; i \in [1,r] $:

\begin{eqnarray} \label{hamderivatives}
\partial_{P^0} H(P,Q)(u)  &=&  \int_{\Sigma_0}  \langle v \circ Q^0(y), u(y) \rangle d\mu_{| \Sigma_0}(y), \\ \nonumber 
\partial_{P^i} H(P,Q)(\Psi) &=&  \int_M \Psi(y) \left( s(Q^i(y)) - \langle v(y) , \nabla Q^i(y) \rangle \right) d\mu(y), \\ \nonumber
\partial_{Q^0} H(P,Q)(u) &=& \int_{\Sigma_0}  \langle [dv]_{|Q^0(y)}(u(y)) , P^0(y) \rangle d\mu_{| \Sigma_0}(y), \\ \nonumber
\partial_{Q^i} H(P,Q)(\Psi) &=& \int_M \Psi(y) [ds]_{|Q^i(y)}P^i(y) - \langle v(y) , P^i(y) \nabla \Psi(y) \rangle d\mu(y).  \nonumber
\end{eqnarray}

\noindent 
Remark that only the last equation really needs to be defined as a distribution and not as a function. Now we can give a sense to the
Hamiltonian equations but only in a weak sense: 

\begin{defi}
An application $\chi \in C^0([0,T],L^1(\Sigma_0,M) \times W^{1,\infty}(M)^r \times L^1(\Sigma_0,\ms{R}^n) \times L^1(M,\ms{R})^r)$
is said to be a weak solution if it verifies for $\Psi \in C_0^{+\infty}([0,T] \times M,\ms{R})$ and $u \in C_0^{+\infty}([0,T] \times M,\ms{R}^n)$:
(we denote $\chi(t)=(Q_t,P_t)$.)
\begin{eqnarray} \label{weak}
\int_0^T \int_M -\partial_t \Psi \; Q_t^i d\mu dt = \int_0^T \partial_{P^i} H(P_t,Q_t)(\Psi) \; dt \; \forall \; i \in [1,r], \\
\int_0^T \int_M -\partial_t \Psi \; P_t^i d\mu dt = -\int_0^T \partial_{Q^i} H(P_t,Q_t)(\Psi) \; dt \; \forall \; i \in [1,r], \\
\int_0^T \int_{\Sigma_0} -\partial_t u \; Q^0_t d\mu_{|\Sigma_0} dt = \int_0^T \partial_{P^0} H(P_t,Q_t)(u) \; dt, \\
\int_0^T \int_{\Sigma_0} -\partial_t u \; P^0_t d\mu_{|\Sigma_0} dt = -\int_0^T \partial_{Q^0} H(P_t,Q_t)(u) \; dt.
\end{eqnarray}
\end{defi}

\subsection{Uniqueness of the weak solutions}

In this paragraph, the uniqueness to the weak Hamiltonian equations is proven, and the proof gives also the
general form of the solutions. This form is closely related to the solution of the variational problem of the 
previous section.

\begin{theo}
Every weak solution is unique and there exists an element of $L^2([0,T],V \times S)$ 
which generates the flow $(\phi_{0,t},\eta_{0,t})$ such that:
\begin{eqnarray} \label{positions}
Q_t^0(x) &=& \phi_{0,t}(x) , x \in \Sigma_0, \\
Q_t^i(u) &=& \eta_{0,t} \circ Q_0^i \circ \phi_{t,0}(u), u \notin \phi_{0,t}(\Sigma_0), i \in [0,n].,
\end{eqnarray}
and for the momentum variables:
\begin{eqnarray} \label{moments}
P_t^0(x) &=& d [\phi_{0,t}]_x^{-1*} (P_0^0(x)) , x \in \Sigma_0, \\
P_t^i(u) &=& P_0^i \circ \phi_{t,0} \Jac (\phi_{t,0}) d[\eta_{t,0}]_{Q_t^i(u)} , u \notin \phi_{0,t}(\Sigma_0).
\end{eqnarray}
\end{theo}

\vspace{0.3cm}

\noindent
\emph{Proof:}
Let $\chi$ a weak solution on $[0,T]$, we introduce 
$$t \mapsto v_t(.)= \frac{1}{\lambda} \int_{\Sigma_0} k(Q^0_t(x),.)P^0_t d\mu_{| \Sigma_0}(x) - \sum_{i=1}^r \int_M  k_V(x,.)P^i_t \nabla Q^i_t d\mu(x),$$
which lies in $L^2([0,T],V)$. This vector field is uniquely determined by the weak solution.
From the preliminaries, we deduce that 
$\phi_{0,t}(x)=\int_0^t v_s(\phi_{s,0}(x)) ds$ is well defined. 
We introduce also,
$$t \mapsto s_t(.)= \sum_{i=1}^r \int_M k_S(Q^i_t,.)P^i_t d\mu(x).$$
For the same reasons, we can integrate the flow:
$\eta_{0,t}(x)=\int_0^t s_r(\eta_{r,0}(x)) dr$ is well defined.
Introducing $\tilde Q_t^i(x)=\eta_{t,0} \circ Q_t^i \circ \phi_{0,t}$ for $i \in [1,r]$, 
we obtain, with $S_t \circ \eta_{t,0}(x)=\partial_t \eta_{t,0}(x)$ and $V_t \circ \phi_{0,t}(x)=\partial_t \phi_{0,t}(x)$:

\begin{eqnarray*}
\int_0^T \int_M -\partial_t \Psi \; \tilde Q_t^i d\mu &=&  \int_0^1  \int_M -\partial_t \Psi \; \eta_{t,0} \circ Q^i_t \circ \phi_{0,t} d\mu dt \\
&=& \int_0^T \int_M -\eta_{t,0} \circ Q^i_t \; [\partial_t \Psi] \circ \phi_{t,0} \Jac (\phi_{t,0})d\mu dt \\ 
&=& \int_0^T \int_M -\eta_{t,0} \circ Q^i_t \; (\partial_t [\Psi \circ \phi_{t,0}] - <\nabla \Psi \circ \phi_{t,0},v_t \circ \phi_{t,0}>)\Jac (\phi_{t,0}) d\mu dt \\ 
&=& \int_0^T \int_M (S_t(\tilde Q^i_t)-<\nabla \tilde Q^i_t,V_t> + d\eta_{t,0}(\dot{Q^i_t} \circ \phi_{0,t})) \Psi \; d\mu dt \\ 
&=& \int_0^T \int_M  (S_t(\tilde Q^i_t)-<\nabla \tilde Q^i_t,V_t> + \\
&\phantom{}& d\eta_{t,0}(-< \nabla Q^i \circ \phi_{0,t} ,v \circ \phi_{0,t} > + s(Q^i \circ \phi_{0,t})) ) \Psi \; d\mu dt.
\end{eqnarray*}

\noindent
The cancelation of the equation above relies on the group relation of flows of vector fields.
We have the equality:
$$ S_t + d\eta_{t,0} (s_t \circ \eta_{0,t}) = 0, $$
then the first and last terms cancel.
The remaining terms cancel too because of the relations:
\begin{eqnarray*}
\nabla \tilde Q^i_t &=& d\phi_{0,t}^* \left( d\eta_{t,0} (\nabla Q_t^i \circ \phi_{t,0}) \right) , \\
v_t &+& d\phi_{0,t} (V_t \circ \phi_{t,0}) = 0.
\end{eqnarray*}
Then, we conclude:

$$\int_0^T \int_M -\partial_t \Psi \; \tilde Q_t^i \; d\mu dt=0.$$

\vspace{0.2cm}
\noindent
Introducing $\Psi(t,x) = \lambda(t) \gamma(x)$, with $\lambda \in C^{+\infty}([0,T])$ and $\gamma \in C^{+\infty}(M)$,
we have:
\\
$\int_0^T - \lambda^{'}(t) (\int_M \gamma(x) \; \tilde Q^i_t d\mu) dt=0$,
hence: $\int_M \gamma(x) \; \tilde Q^i_t d\mu = \int_M \gamma(x) \; \tilde Q^i_t d\mu$,
i.e.
$\tilde Q_t^i = \tilde Q_0^i = Q_0^i$, and:

$$ Q_t^i = \eta_{0,t} \circ Q_0^i \circ \phi_{t,0}. $$

\vspace{0.3cm}
\noindent
Now, we introduce for $i \in [1,r]$, $\tilde P_t^i(.) = \frac{P_t^i \circ \phi_{0,t} \Jac (\phi_{0,t}(.))}{ d[\eta_{t,0}]_{\eta_{0,t} \circ Q_0^i(.)}}$,
this quantity is well defined because of the inversibility of the flow of $s_t$ and $v_t$.
Remark that $\frac{\Jac (\phi_{0,t}(.))}{ d[\eta_{t,0}]_{\eta_{0,t} \circ Q_0^i(.)}}$
is differentiable almost everywhere because $Q_0^i$ is \lip on $M$. We want to prove that
$\int_0^1 \int_M -\partial_t \Psi \; \tilde P_t^i d\mu dt=0$ with $\Psi \in C^{+\infty}_0(M)$, which leads to:
$\tilde P_t^i = \tilde P_0^i = P_0^i$, and then we are done. 

\vspace{0.2cm}

\noindent
To prove the result, we first use the change of variable $y=\phi_{0,t}(x)$, this is a straightforward calculation,
we will also use the equality: $d[\eta_{t,0}]_{\eta_{0,t}(.)}  d\eta_{0,t} (.) = 1$.
\begin{eqnarray}
\int_0^T \int_M \Psi \partial_t \tilde P_t^i \; d\mu dt &=& -\int_0^T \int_M \tilde P_t^i \partial_t \Psi \circ \phi_{t,0} \Jac(\phi_{t,0})\; d\mu dt, \nonumber \\ 
\int_0^T \int_M \Psi \partial_t \tilde P_t^i \; d\mu dt &=& -\int_0^T \int_M \frac{P_t^i}{d[\eta_{t,0}]_{Q_t^i(.)}} \partial_t \Psi \circ \phi_{t,0}\; d\mu dt, \nonumber \\ 
\int_0^T \int_M \Psi \partial_t \tilde P_t^i \; d\mu dt &=& \int_0^T \int_M -P_t^i \partial_t (\frac{\Psi \circ \phi_{t,0}}{d[\eta_{t,0}]_{Q_t^i(.)}}) + \frac{P_t^i}{d[\eta_{t,0}]_{Q_t^i(.)}} <\nabla \Psi _{|\phi_{t,0}},-d\phi_{t,0}(v_t)> \nonumber \\ 
 &+& P_t^i \Psi \circ \phi_{t,0}\; \partial_t (d[\eta_{0,t}]_{Q_0^i \circ \phi_{0,t}}) d\mu dt, 
\end{eqnarray}

\noindent
The third term of the last equation can be rewritten:
\begin{eqnarray} \label{intermed}
&\phantom{1}& \int_0^T \int_M P_t^i \Psi \circ \phi_{t,0}\; \partial_t (d[\eta_{0,t}]_{Q_0^i \circ \phi_{0,t}}) \; d\mu dt = \nonumber \\
&\phantom{1}& \int_0^T \int_M P_t^i \Psi \circ \phi_{t,0}\; (<\nabla (d[\eta_{0,t}]_{Q_0^i}),-d\phi_{t,0}(v_t)> + d[s_t]_{Q_t^i}d[\eta_{0,t}]_{Q_0^i \circ \phi_{0,t}}) \; d\mu dt.
\end{eqnarray}

\noindent
Now, we can apply the hypothesis on $\partial P_t^i$, the first term of the expression is equal to:

\begin{eqnarray}
\int_0^T \int_M -P_t^i \partial_t (\frac{\Psi \circ \phi_{t,0}}{d[\eta_{t,0}]_{Q_t^i(.)}}) \; d\mu dt &=& \int_0^T \int_M P_t^i <\nabla (\Psi \circ \phi_{t,0}d[\eta_{0,t}]_{Q_0^i \circ \phi_{0,t}}),v_t>  \nonumber \\ 
&-& \Psi \circ \phi_{t,0} \; d[\eta_{0,t}]_{Q_0^i \circ \phi_{0,t}} d[s_t]_{Q_t^i} P_t^i \; d\mu dt, \\ 
\int_0^1 \int_M P_t^i <\nabla (\Psi \circ \phi_{t,0}d[\eta_{0,t}]_{Q_0^i \circ \phi_{0,t}}),v_t>  \; d\mu dt &=& \int_0^1 \int_M P_t^i <d\phi_{t,0}^*(\nabla \Psi_{|\phi_{t,0}}),v_t> d[\eta_{0,t}]_{Q_0^i \circ \phi_{0,t}} \nonumber \\ 
&+& P_t^i <d\phi_{0,t}^* \nabla (d[\eta_{0,t}]_{Q_0^i}),v_t>\; d\mu dt.
\end{eqnarray}

\noindent
All the terms of the equation cancel together, so we obtain the result.

\vspace{0.3cm}

\noindent
With $\tilde P_t^0(x) = d [\phi_{0,t}]_x^* P_t^0(x)$ and $\tilde Q_t^0 = \phi_{t,0} \circ Q_0^0$,
we get the result for the last two terms of the system in the same way than the preceding equations, but it is even easier.
$\square$

\vspace{0.3cm}
\noindent
Remark that we only have to suppose the weak solution is $L^2$ to obtain the result. More than uniqueness, 
we know that the weak solution "looks like" a variational solution of our initial problem.

\subsection{On the existence of weak solutions}

We consider in this section a Hamiltonian equation which includes our initial case, for which we have proved 
existence results. Namely, we can rewrite our result on the last section in terms of existence of weak
solution of the system (\ref{ham}).

\begin{prop}
Let $(U_1,\ldots,U_n)$ be a partition in \lip domain of $M$, $\Sigma_0 = \cup_{i=1}^{n} \partial U_i$.
For any intial data,  $I_0 \in W^{1,\infty}(M)^r$, $Q_0 = (\Sigma_0,I_0)$
and $P_0 \in L^{\infty}(\Sigma_0,\ms{R}^n) \times L^1(M,\ms{R})^r$ such that $\text{Supp}(P_0^i) \subset U_i$,
then there exists a solution to the Hamiltonian equations. 
\end{prop}

\noindent
Remark that this solution has the same structure than a variational solution of our initial problem.
In this case, all the momenta ($i \in [1,r]$) can be viewed as one momentum: with $P_t = \sum_{i=1}^r P_t^i$,
we have all the information for the evolution of the system. 
\\
Now, we can say a little bit more on the general Hamiltonian equations. We will not give a proof here
of the existence if we relax the condition on the support of $P^i_0$, but the reader can convince himself that
the existence is somehow a by-product of the last section. 
\\
Summing up our work at this point, from a precise variational problem we obtain generalized Hamiltonian equations,
for which we can prove results on existence and uniqueness. A natural question arises then, from what type of variational
problems could appear these solutions? The answer could be based on the remark: because the decomposition we choose is the direct product of spaces,
we can put a sort of product metric on it.
A simple generalized minization is obtained by modifying the equation \eqref{minimisation}:

\begin{equation} \label{minimisation2}
E(\eta,\phi)=D(Id,(\eta,\phi))^2 + \sum_{i =1}^r \frac{1}{\sigma_i^2} \| \eta \circ (I_0^i \mt{1}_{U_i}) \circ \phi^{-1} -I^i_{targ}\|_{L^2}^2 ,
\end{equation}
with for $i \in [1,r]$ $I_0^i \in W^{1,\infty}(M)$, and $I_0^0 \in \text{Im}(M)$.

\section{Conclusion}
The main point of this paper is the derivation lemma which
may be of useful applications. This technical lemma gives a larger framework
to develop the \emph{large deformation diffeomorphisms} theory.
The action on the level lines is far to be none of interest but
we aim to obtain numerical implementations of the contrast term applied to smooth
images. 
Finally, the interpretation as a Hamiltonian system through optimal control theory ends up with giving a proper 
understanding of the momentum map.
To go further, the technical lemma seems to be easily enlarged to rectifiables domains, 
and there may be a useful generalization to $SBV$ functions. This would enable a generalization
of a part of this work to $SBV$ functions. But to understand the weak Hamiltonian formulation
would have been much more difficult within the $SBV$ framework. From the numerical point of view,
some algorithms that are currently developed to treat the evolution of curves
could be used efficiently but they need strong developments.

\section{Proof of the lemma} \label{pol}
After recalling some classical facts about \lip functions, we prove the derivation lemma:

\begin{lem} 
Let $U,V$ two bounded \lip domains of $\mt{R}^n$. 
Let $X$ a \lip vector field on $\ms{R}^n$ and $\phi_t$ the associated flow. 
Finally, let $g$ and $f$ \lip real functions on $\mt{R}^n$.
Consider the following quantity depending on $t$,
$$ J_t = \int_{\phi_t(U)} f \circ \phi_{t}^{-1} g \mt{1}_V d\mu,$$
where $d\mu$ is the Lebesgue measure, then
\begin{equation} 
\partial_{t|t=0^+} J_t = \int_{U} -<\nabla f,X> g \mt{1}_V d\mu + \int_{\partial U} <X,n> fg  \tilde{\mt{1}}_V(X) d\mu_{|\partial U}.
\end{equation}
with $ \tilde{\mt{1}}_V(X)(y)= \lim_{\epsilon \mapsto 0^+} \mt{1}_{\bar V}(y+\epsilon X) $, if the limit exists, $0$ elsewhere.
And we denote by $d\mu_{|\partial U}$ the measure on $ \partial U$ and $n$ the outer unit normal of $\partial U$.
\end{lem}

\noindent
We will use,
\begin{theo}\label{Rademacher}
\textbf{Rademacher's theorem}
\\
Let $f: U \mapsto \ms{R}^n$ a \lip function defined on an open set $U \subset \ms{R}^n$, 
then $f$ is differentiable $\mu$ a.e.
\end{theo}

\begin{theo} \label{implicitfunctions}
Let $F: \ms{R}^d \times \ms{R}^n \mapsto \ms{R}^d$, \lip continuous on a neighborhood of $(v_0, y_0)$ and $F(v_0,y_0)=0$.
Suppose that $\partial_1 F(v_0,y_0)$ exists and is invertible. Then there exists a neighborhood $W$ of  $(v_0, y_0)$, on which there
exists a function $g: \ms{R}^n \mapsto \ms{R}^d$ such that in $W$:

\begin{itemize}
\item $g(y_0)=v_0$.
\item $F(g(y),y)=0 \; for \; (g(y),y) \in W$.
\item $|g(y)-g(y_0)| \leq c |y-y_0|$, 
\end{itemize}
with $c=1+(\text{Lip}(F)+1)\|\partial_v F(v_0,y_0)^{-1} \|$.
\end{theo}

\noindent
This theorem can be found in \cite{1022.49018} and in a more general exposition than we will use hereafter.

\vspace{0.3cm}

\noindent
Now an obvious lemma of derivation under the integral,

\begin{lm}
Let $f$ a \lip function defined on an open set $U \subset \ms{R}^n$.
Let $X$ a \lip vector field on $\ms{R}^n$ of compact support and $\phi_t$ the associated flow. 
Consider $J_t = \int_U f \circ \phi_t dx$, then:
\begin{equation} 
\partial_{t|t=0} J_t = \int_{U} <\nabla f,X> dx.
\end{equation}
\end{lm}

\noindent
\emph{Proof:}
Using Rademacher's theorem, this is a staightforward application of dominated convergence theorem. 
Note that under the condition that $f$ is Lipschitz, if both $f$ and $\nabla f$ are integrable and $X$ is a bounded vector field,
we can relax the hypothesis of a compact support for the vector field, which is replaced here by the integrability condition. 
$\square$

\vspace{0.2cm}
\noindent
We will need the following characterization of derivation for real functions to prove the lemma \ref{implicite} .

\begin{lm} \label{sursous2}
Let $w: \ms{R}^n \mapsto \ms{R}$ a function,
then $w$ is differentiable in $x \in \ms{R}^n$ if and only if
there exist $f$ and $g$ two $C^1$ functions and a neighborhood $V$ of $x$,
such that 
$f(x)=g(x)$  and if $y \in V$,
\begin{equation}
g(y) \leq w(y) \leq f(y).
\end{equation}
\end{lm}

\noindent
\emph{Proof:}
Suppose $w$ differentiable in $x$,
it suffices to prove that there exists $f$ $C^1$ such that $w \leq f$ 
in a neighborhood of $x$. (To obtain $g$, consider then $-w$.)  
we can suppose $x=0$, $w(0)=0$ and $w'(0)=0$
then $\lim_{y \mapsto x} \frac{w(y)}{|y|} = 0$.
Hence there exists a continuous function $v$ defined on $B(0,r>0)$ such that:
$\frac{w(y)}{|y|} \leq v(y)$ and $v(0)=0$. 
\noindent
With the notation $|x|$ for the euclidean norm in $\ms{R}^n$,
let $m(|x|)=\sup_{|y| \leq |x|} v(y)$, we have $m(|x|) \geq v(x)$ for $x \in B(0,r)$.
The function $m$ is non decreasing and continuous with $m(0)=0$.
At last, let $f(x)=\int_{|x|}^{2|x|} m(t) dt$, then $w(x) \leq |x|m(|x|) \leq f(x)$.
Moreover, the fact $f$ is $C^1$ is straightforward to verify.
\\
Suppose $f$ and $g$ are $C^1$, and denote by $f'$ and $g'$ their derivative in $x$, then
$f'(x)=g'(x)$ since $f-g \geq 0$ and has a minimum in $x$.
On $V$ we have: $$g(x+h)-g'(x).h \leq w(x+h)-g'(x).h \leq f(x+h)-f'(x).h.$$ 
Hence, $w(x+h)-g'(x).h=o(h)$ and $w$ is differentiable. (Remark that we only use the fact
that $f$ and $g$ are differentiable in $x$.)
$\square$

\vspace{0.2cm}
\noindent
Using the lemma above, we study the deformation of an epigraph of a \lip function under the 
action of a vector field, which leads to study the deformation of the graph of the function:

\begin{lm} \label{implicite}
Let $\phi_t$ the flow of the vector field \lip $X$ on $\ms{R}^n$ (with $\|X\|$ bounded on $\ms{R}^n$). Let $V=\{(x,z) \in \ms{R}^{n-1} \times \ms{R} | z>w(x)\}$ with 
$w$ a \lip function, and $w_t(x)=\inf \{ z | (x,z) \in  \phi_t(V) \}$.
\\
Then, a.e.
$$\partial_{t | t=0} w_t(x) = -< \nabla w(x) , p_1(X(x,w(x)))> + p_2(X(x,w(x))),$$
with $p_1$ and $p_2$ orthogonal projections respectively on $\ms{R}^{n-1} \times 0$ and $0^{n-1} \times \ms{R}$.
\end{lm}

\noindent
\emph{Proof:}
Remark that $w_t$ is well defined for all $t$ by connexity reason, but $w_t$ might be discontinuous for $t$ large enough.
However it is \lip continuous for t in a neighborhood of $0$: 
\noindent
we first apply the implicit function theorem for \lip maps to the function:
$$F(x,t)=p_1(\phi_t(x,w(x)))-x_0,$$
note that $\partial_x F(x_0,0) = Id_{| \ms{R}^{n-1}}$, so we obtain for each $x_0 \in \ms{R}^{n-1}$
a function $x_0 :t \mapsto x_0(t)$ such that $x_0$ is \lip and the equation $F(x,t)=0 \Leftrightarrow x=x_0(t)$
on a neighborhood of $(x_0,0)$. Note that the implicit function theorem in \cite{1022.49018} gives only existence but not 
uniqueness. We develop now the uniqueness.
The \lip condition on $w$ can be written with the cone property. Let $(A,B)$ two points on the graph of $w$, then 
$|y_A-y_B| < M |x_A-x_B|$, for a \lip constant. This open condition is then verified in a neighborhood of $\{t=0\}$.
We see that, $F(x_1,t)=F(x_2,t)$ implies 
$ \phi_t(x_1,w(x_1))=\phi_t(x_2,w(x_2))$, hence $x_1=x_2$.

\vspace{0.3cm}
\noindent
If $w$ is $C^1$, we get by implicit function theorem the first derivative of $x_0(t)$: 
\begin{eqnarray*}
\partial_{t|t=0} x_t&=& -p_1(X(x_t,w(x_t))).
\end{eqnarray*}
We deduce,
\begin{equation} \label{impl}
w_t(x)=p_2(\phi_t(x_t,w(x_t))),
\end{equation}
and that $t \mapsto w_t(x_0)$ is a \lip function for each $x_0$. 
In the $C^1$ case, we get by differentiation of the equation \eqref{impl},
$$\partial_{t_{|t=0}}w_t(x_0)=-< \nabla w(x) , p_1(X(x,w(x)))> + p_2(X(x,w(x))).$$
\\
We now observe that there is an obvious monotonicity in $w$ of $w_t$.
Indeed, if $w \leq v$ then $w_t \leq v_t$.
We then use the lemma \ref{sursous2} to prove the result in the case 
$w$ is \lip .
\noindent
Let $x$ such that $w$ is upper and lower approximated by $C^1$ functions:
let $u$ and $v$ such that 
$u(x)=w(x)=v(x)$, and $u \leq w \leq v$.
We obtain:
\begin{equation}
\frac{1}{t}(u_t(x)-u(x)) \leq \frac{1}{t}(w_t(x)-w(x)) \leq \frac{1}{t}(v_t(x)-v(x)),
\end{equation}
We deduce the result:
$$\lim_{t \mapsto 0} \frac{1}{t}(w_t(x)-w(x)) = -< \nabla w(x) , p_1(X(x,w(x)))> + p_2(X(x,w(x))),$$
for all the points of derivability of $w$, i.e. almost everywhere since $w$ is Lipschitz. 
$\square$

\vspace{0.2cm}
\noindent
Now, we prove the following lemma, which can be seen as a consequence of the coarea formula.
It will be used in the proposition \ref{lemmedebase}.

\begin{lm} \label{niv}
Let $w: \ms{R}^n \mapsto \ms{R}$, an a.e. differentiable function, 
$A=w^{-1}(\{0\})$ and $B= \{x | \nabla w(x) \neq 0 \}$ then $\mu(A \cap B)=0$.
\end{lm}

\noindent
\emph{Proof:}
For $n=1$, the lemma is obvious because the point of $A \cap B$ are isolated with Taylor formula.
For $n>1$, we generalize with Fubini's theorem:
$$ A \cap B = \bigcup_{i=1}^n A \cap B_i,$$
with $B_i = \{ x \in \ms{R}^n | <\nabla f(x),e_i> \neq 0\}$.
To prove the result, it suffices to see that $\mu(A \cap B_n) = 0$.
Consider $(x,t) \in \ms{R}^{n-1} \times \ms{R}$, $f_x(t)=f(x+t)$ and $D_x=x \times \ms{R}$. 
Then we apply the case $n=1$ to the function $f_x$:
$\mu(A \cap B_n \cap D_x)=0$, and with Fubini's theorem, $\mu(A \cap B_n) = 0$. 
$\square$

\vspace{0.2cm}
\noindent
Remark that this lemma can be applied to a \lip function.
Below lies the fundamental step to prove the derivation lemma.

\begin{prop} \label{lemmedebase}
Let $w: \mt{R}^{n-1} \mapsto \ms{R}$ a \lip function. Let
$V:= \{(x,y)|y>w(x)\}$ and $U:=\{(x,y)| y>0 \}$.
Let $X$ a \lip vector field on $\ms{R}^n$ and $\phi_t$ the associated flow. 
Finally, let $g$ and $f$ \lip real functions on $\mt{R}^n$ of compact support.
Consider the following quantity depending on $t$,
$$ J_t = \int_{\phi_t(U)} f \circ \phi_{t}^{-1} g \mt{1}_V d\mu,$$
where $d\mu$ is the Lebesgue measure also denoted by $dx$, then
\begin{equation} 
\partial_{t|t=0^+} J_t = \int_{U} -<\nabla f,X> g \mt{1}_V dx + \int_{\partial U} <X,n> fg  \tilde{\mt{1}}_V(X) d\mu_{|\partial U}.
\end{equation}
with $ \tilde{\mt{1}}_V(X)(y)= \lim_{\epsilon \mapsto 0^+} \mt{1}_{\bar V}(y+\epsilon X) $, if the limit exists, $0$ elsewhere.
And we denote by $d\mu_{|\partial U}$ the measure on $ \partial U$ and $n$ the outer unit normal.
\end{prop}

\begin{figure}[htbp]
	\centering
		\includegraphics{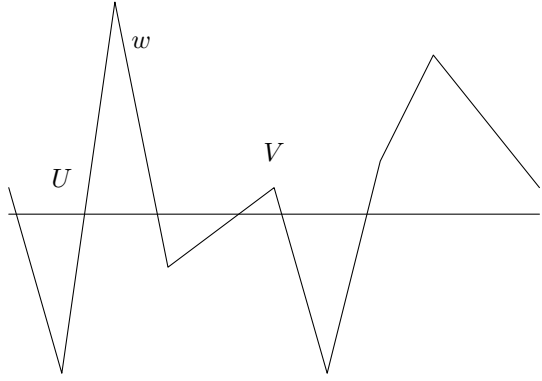}
	\caption{\footnotesize The main case}
	\label{fig:premiercas}
\end{figure}


\noindent
\emph{Proof:}
The first case to treat is when $w>0$, we can then integrate on $V$ instead of $\phi_t(U)$.
$$J_t = \int_{V} f \circ \phi_{t}^{-1} g \mt{1}_V d\mu,$$ we differentiate under the integral, we get:
$$ \partial_{t|t=0} J_t = \int_{U} -<\nabla f,X> g \mt{1}_V d\mu ,$$
this is the formula because the second term is null. In the following, we have to use this case.

\vspace{0.3cm}

\noindent
To treat the general case, we first do a change of variable:
$$J_t=\int_{\ms{R}^{n-1}}\int_{w_t(x)^+}^{+\infty}f g \circ \phi_t \Jac (\phi_t) d\mu.$$
\noindent
We introduce some notations:
\begin{eqnarray*}
x^+ &=& \max (0,x) =H(x),\\
\nabla H(x)(v) &=& 0, x<0,\; or \;x=0 \;et\; v<0, \\
\nabla H(x)(v) &=& v \; ,elsewhere, \\  
w_t(x)&=&\inf \{ z | (x,z) \in  \phi_{-t}(V) \}.
\end{eqnarray*}

\noindent
Let $p_1$ and $p_2$ the orthogonal projections respectively on $\ms{R}^{n-1} \times \{0\}$ and $0^{n-1} \times \ms{R}$.
With the lemma \ref{implicite}, 
$$\partial_{t | t=0^+} w_t(x) = < \nabla w(x) , p_1(X(x,w(x)))> - p_2(X(x,w(x))).$$
Consequently,
$$\partial_{t | t=0^+} w_t(x)^+ = \nabla H(w(x))(< \nabla w(x) , p_1(X(x,w(x)))> - p_2(X(x,w(x)))).$$
\noindent
Using $f(y) (<\nabla g,X> + \div(X)g)=\div(fgX)-g<\nabla f,X>$, we get:
\begin{eqnarray*}
\partial_{t|t=0^+} J_t = \int_{\ms{R}^{n-1}} -\partial_{t | t=0^+} w_t(x)^+ f(x,w(x)^+) g(x,w(x)^+) dx \\
- \int_{\ms{R}^{n-1}} \int_{w_t(x)^+}^{+\infty} <\nabla f(x,z),X(x,z)> g(x,z) dxdz\\
+\int_{\partial (U \cap V)}  f g  <X,n> d\mu_{|\partial (U \cap V)}.
\end{eqnarray*}
Here $n$ is the outer unit normal of $\partial(U \cap V)$. 
Rewrite the last term:
\begin{eqnarray*}
\int_{\partial (U \cap V)}  f g  <X,n> d\mu_{|\partial (U \cap V)}&=&\int_{\partial V \cap U}  f  g  <X,n> d\mu_{|\partial V} \\
&+&\int_{w^{-1}(]-\infty,0])}f(x,0)g(x,0)<X,n> dx.
\end{eqnarray*}
\noindent
In a neighborhood $C$ of $x$ such that $w(x) > 0$, we have 
demonstrated that the formula holds, so the first term in the equation above is equal to:
$$ -\int_{\partial V \cap U}  f  g  <X,n> d\mu_{|\partial V}-\int_{w^{-1}(]-\infty,0])} \partial_{t | t=0} w_t(x)^+ f(x,0) g(x,0) dx.$$

\noindent
Moreover, on the set $F=\{ x: w(x) = 0\}$, we have, with the lemma \ref{niv}, 
a.e. $\nabla w =0$. Then, we have: $\partial_{t | t=0^+} w_t(x)^+ = <X,n>^+$ , a.e.
On the set $G=\{ x: w(x) < 0\}$, we have: $\partial_{t | t=0^+} w_t(x)^+ = 0$.
We now get the result with:
\begin{eqnarray*}
\int_{w^{-1}(]-\infty,0])} (<X,n>-\partial_{t | t=0} w_t(x)^+) f(x,0) g(x,0) dx=  \\
\int_{w^{-1}(]-\infty,0])} \tilde{\mt{1}}_{V}(X)(x,0) <X,n> f(x,0) g(x,0) dx.
\end{eqnarray*}
Indeed, if $<X, n> \neq 0$ the result is straightforward because the limit exists in the definition of $\tilde{\mt{1}}_V$. 
If $<X, n> = 0$, the contribution is null. 
$\square$

\vspace{0.2cm}
\noindent
Our goal is to prove the formula for \lip open sets, we present some definitions. 

\begin{defi} \label{lip}
$\;$
\\
An open set $U \neq \emptyset$ of $\ms{R}^n$ is said to be locally \lip  if
for each $x \in \partial U$, there exist:
\begin{itemize}
 \item an affine isometry $I$, of $\ms{R}^n$, 
 \item an open neighborhood $V(x)$ of $x$,
 \item a \lip function $w$ defined on $\ms{R}^{n-1}$ with a \lip constant $K(x)$
\end{itemize}
such that, $$I(V(x) \cap U) = I(V(x)) \cap \{(x,y) \in \ms{R}^{n-1} \times \ms{R} |  y>w(x)\}$$
\\
If the constant $K(x)$ can be chosen independent of $x$, $U$ is said to be Lipschitz. 
\end{defi}

\begin{rem}
$\;$
\begin{enumerate}
\item An open bounded set of $\ms{R}^n$ which is locally \lip is also \lip . 
\item By Rademacher's theorem, the outer unit normal $n(x)$ exists for $\mathcal{H}^{n-1}$ a.e. $x \in \partial U$.
\item We will say that $(V(x),I)$ trivializes the \lip domain in $x$.
\end{enumerate}
\end{rem}

\vspace{0.2cm}
\noindent
The three lemmas below prove that one can describe a \lip domain
in many systems of coordinates. This is a key point to understand the 
two boundaries at a point of intersection and enables to use the proposition
\ref{lemmedebase}.

\begin{lm}
Let $\psi$ a $C^1$ diffeomorphism of $\ms{R}^n$ and $V$ a \lip domain,
then $\psi(V)$ is a \lip domain.
\end{lm}

\noindent
\emph{Proof:}
The proof is straightforward with the characterization of \lip domains
with the uniform cone property, which can be found in \cite{1002.49029}.
$\square$

\begin{lm} \label{lm0}
Let $(e_1,\ldots,e_n)$ an orthonormal basis of $\ms{R}^n$, and $w$ a \lip function defined on $\ms{R}^{n-1}$ 
of \lip constant $M$.
Let $U$ the \lip open set which is above the graph of $w$:
$$U:= \{(x,y) \in \ms{R}^{n-1} \times \ms{R} |  y>w(x)\},$$
then for each $n$ in the open cone $C:= \{ n=(x,y) \in \ms{R}^{n-1} \times \ms{R_+} |y>M |x| \}$
one can trivialize the boundary of $U$ through the graph of a function defined on $n^{\perp}$ (with an orthonormal basis).
Moreover, this function is \lip . 
\end{lm}


\begin{figure}[htbp]
	\centering
		\includegraphics{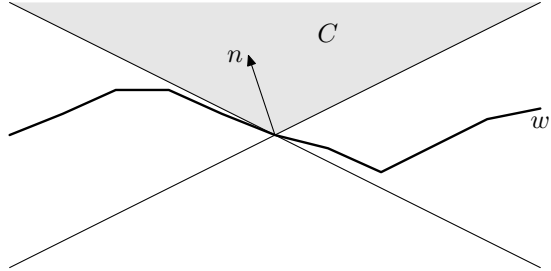}
	\caption{\footnotesize Trivializing with respect to the orthogonal hyperplan to $n$.}
	\label{fig:trivialisation}
\end{figure}

\noindent
\emph{Proof:}
Note that, $n$ can be represented with the angle between the hyperplan and the vector orthogonal to this hyperplan, and also
the Lipschitz constant can be represented as the tangente of such an angle. 
Let two points $a,b \in \ms{R}^n$ which belong to the $\partial U$,
then $b-a$ and $n$ are not colinear, because of the \lip property.
As a consequence, $\partial U$ is defined as the graph of a function $\tilde w$ on $n^{\perp}$.
And, one can verify that, if $n$ is normalized, a \lip constant for $\tilde w$ is equal to: $tan(|\theta_1-\theta_2|)$,
if $\theta_1$ is the angle of $n$ and $\theta_2$ associated to the \lip constant.   
$\square$

\begin{lm} \label{lm1}
Let $U \subset \ms{R}^n$ a \lip domain with $0 \in \partial U$, then there exist a neighborhood $V$ of $0$,
$w$ a \lip function defined on $\ms{R}^{n-1}$ and a linear transformation $A$ such that:
$(e_1,\ldots,e_{n-1}) \subset Ker(A-Id)$ and 
$$V \cap A(U) = V \cap \{(x,y) \in \ms{R}^n \times \ms{R}| y > w(x)\}.$$
\end{lm}

\noindent
Here is the illustration of the idea driving the proof.

\begin{figure}[htbp]
	\centering
		\includegraphics{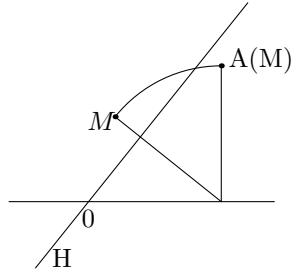}
	\caption{\footnotesize The linear transformation $A$}
	\label{fig:aff}
\end{figure}

\noindent
\emph{Proof:}
In some coordinates, we write $U$ as the epigraph $G^+$ of a \lip function $v$ on a hyperplan $H$ in a neighborhood $W$ of $0$.
We note by $n$ a normal vector to $H$. We face two cases:

\vspace{0.3cm}

\noindent
$\bullet$ If $<n,e_n> \neq 0$, let $A \in L(\ms{R}^n)$ defined by: 
$A(n)=e_n$ and $(e_1,\ldots,e_{n-1}) = Ker(A-Id)$. 
Denote by $p_2$ the orthogonal projection on $(e_1,\ldots,e_{n-1})$. Let $M \in G$ defined by its projections $z$ on $H$:
$M=z+w(z)n$, by definition of $A$, $A(M)= (p_2(z),w(z))+A(z-p_2(z))$. But $z-p_2(z)=\lambda(z) n$ with $\lambda \in (\ms{ R}^{n})'$ 
so we obtain:
$A(M)= (p_2(z),w(z) + \lambda)$.
Also $p_{2|H}$ is a linear isomorphism, we note the inverse $p_2^{-1}$,
then with the change of variable $x=p_2(z)$, we get $$G=\{(x,\lambda(p_2^{-1}(x)) + w \circ p_2^{-1}(x)) | x \in p_2(H \cap W)\}.$$  
$w \circ p_2^{-1}$ is clearly \lip and we obtain the lemma in this case.

\vspace{0.3cm}

\noindent
$\bullet$ If $<n,e_n> = 0$, we can choose by lemma \ref{lm0} another system of coordinates for which 
we fall in the first case, and the lemma is demonstrated.
$\square$

\vspace{0.2cm}
\noindent
We present a smooth ($C^1$) version of the derivation lemma.

\begin{prop} \label{lisse}
Let $U$ a bounded $C^1$ domain of $\mt{R}^n$ and $V$ a bounded \lip domain. 
Let $X$ a \lip vector field on $\ms{R}^n$ and $\phi_t$ the associated flow. 
Finally, let $g$ and $f$ \lip real functions on $\mt{R}^n$.
Consider the following quantity depending on $t$,
$$ J_t = \int_{\phi_t(U)} f \circ \phi_{t}^{-1} g \mt{1}_V d\mu,$$
where $d\mu$ is the Lebesgue measure also denoted by $dx$, then
\begin{equation}
\partial_{t|t=0^+} J_t = \int_{U} -<\nabla f,X> g \mt{1}_V dx + \int_{\partial U} <X,n> fg  \tilde{\mt{1}}_V(X) d\mu_{|\partial U}.
\end{equation}
with $ \tilde{\mt{1}}_V(X)(y)= \lim_{\epsilon \mapsto 0} \mt{1}_{\bar V}(y+\epsilon X) $, if the limit exists, $0$ elsewhere.
And we denote by $d\mu_{|\partial U}$ the measure on $ \partial U$ and $n$ the outer unit normal.
\end{prop}

\vspace{0.3cm}
\noindent
\emph{Proof:}
Let $K$ a \lip constant of $V$. 
Applying the definition of a $C^1$ domain with the compact boundary of $U$, 
there exist a finite covering $W_1,\ldots,W_k$ of $\partial U$ with open balls and $(\Psi_1,\ldots,\Psi_k)$ $C^1$ diffeomorphisms such that
for each $i \in [1,k]$,
$$ \Psi _i (W_i \cap U) = \{ (x,y) \in  \ms{R}^{n-1} \times \ms{R} |  y>0 \} \cap \Psi _i (W_i),$$  
also one has: 
$$\Psi_i (W_i \cap \partial U) = \ms{R}^{n-1} \times \{0\} \cap \Psi_i (W_i).$$  
\noindent
Let $(\theta_0,\ldots,\theta_n,\theta_{k+1})$ a partition of unity associated to the family
$$(W_0=\ms{R}^n \setminus \bar U,W_1,\ldots,W_k,W_{k+1}=U).$$ 
It means:
\begin{itemize}
\item $0 \leq \theta_i \leq 1, \forall i \in [0,k+1]$ and $\sum_{i \in [0,k+1]} \theta_i = 1$ on $\ms{R}^n$. 
\item $Supp \; \theta_i \subset W_i$ for $i \in [1,k+1]$.
\item $Supp \; \theta_0 \subset \ms{R}^n \setminus \bar U$.
\end{itemize}

\vspace{0.3cm}

\noindent
Through the change of variable $y=\phi_t(x)$, the quantity is:
$$J_t = \sum_{i=1}^{k+1} \int_{U} f \theta_i g \circ \phi_t \mt{1}_{V} \circ \phi_t \Jac (\phi_t) d\mu. $$

\noindent
Four cases appear: 
\begin{itemize}
\item $\bar U \cap W_i \subset V$
\item $W_i \cap U \subset \ms{R}^n \setminus V$
\item $\bar V \cap W_i \subset U$
\item $\partial V \cap \partial U \cap W_i \neq \emptyset$.
\end{itemize}

\vspace{0.3cm}

\noindent 
In the first case, the formula is the result of the derivation under the integral, which is allowed because 
$g$ is Lipschitz. 
$$ \partial_{t|t=0^+} J_t = \int_{W_i} f(y) (<\nabla g,X> + \div(X)g)d\mu, $$
with, $f(y) (<\nabla g,X> + \div(X)g)=\div(fgX)-g<\nabla f,X>$, and applying Stokes theorem
true for a rectifiable open set and \lip functions, we obtain the result. 
\\
In the second case, the quantity is null for $t$ sufficiently small. So the formula is obvious.
\\
In the third case, we can integrate on $V$ instead of $\phi_t(U \cap W_i)$:

$$J_t = \int_{V} f \circ \phi_{t}^{-1} g \mt{1}_V d\mu,$$ we differentiate under the integral, we get:
$$ \partial_{t|t=0^+} J_t = \int_{U} -<\nabla f,X> g \mt{1}_V d\mu ,$$
because the second term of the formula is null.
\\
We deal hereafter with the last case:
as $\Psi_i$ is a $C^1$ diffeomorphism, $\Psi_i(V)$ is also Lipschitz. 
Consequently, we can find a finite covering $B_1, \ldots , B_m$ of $W_i$, for which one of the following conditions holds:

\begin{itemize} 
\item $B_i \subset V$
\item $B_i \subset \ms{R}^n \setminus \bar V$
\item $B_i \cap \partial V \neq \emptyset$ and there exists $I$ such as $(\Psi_i(B_i),I)$ trivializes the \lip domain $\Psi_i(V)$.
\end{itemize}

\noindent
In the first two cases, we have already demonstrated that the formula is true. 
\\
With the lemma \ref{lm1}, we know that after a linear transformation $A$ which is the identity on $\ms{R}^{n-1} \times \{ 0\} $, 
the \lip domain can be represented as the epigraph of a \lip function defined on $\ms{R}^{n-1}$. We replace
$\Psi_i$ by $A \circ \Psi_i=\Psi$. We then have the following situations:
$$\Psi(W_i \cap V) = \{ (z,t) | t > w(z) \},$$ 
or $$\Psi(W_i \cap V) = \{ (z,t) | t < w(z) \},$$ 
with $w : \ms{R}^{n-1} \cap B(0,\rho) \mapsto \ms{R}$ a \lip function.
This situation (or the symetric situation which is essentially the same) is treated in the proposition \ref{lemmedebase}. 
\noindent
$\square$

\vspace{0.3cm}
\noindent
We generalize the proposition \ref{lisse} to the case of \lip domains, we need some additional results of approximation:

\begin{theo}\label{approx}
\textbf{$C^1$ approximation}
\\
Let $f: \ms{R}^n \mapsto \ms{R}$ a \lip function.
Then for each $\epsilon > 0$, there exists a $C^1$ function $\bar f : \ms{R}^n \mapsto \ms{R}$ such that:
$$\mu(\{x | \bar f (x) \neq f(x) or D \bar f (x) \neq Df(x) \}) \leq \epsilon .$$
In addition, 
$$\sup_{\ms{R}^n} | D \bar f | \leq C \text{Lip}(f),$$
for some constant $C$ depending only on $n$.
\end{theo}

\noindent
See the proof of \cite{0804.28001}. 

\begin{rem}
A direct consequence of the theorem is that we have, 
$$ \| f-\bar f \|_{\infty} \leq 2 \max(C,1) \text{Lip}(f) \sqrt{n-1} \epsilon^{\frac{1}{n-1}}.$$
On each cube of volume $\epsilon$ there exists a point where the two functions are equal, 
then we deduce easily the claimed bound. Thus, we get also
$$ \mu(\{(x,y)| \bar f (x)<y< f(x) \; or \; f(x) < y < \bar f(x) \}) \leq 2 \max(C,1) \; \text{Lip}(f) \sqrt{n-1} \epsilon^{1 + \frac{1}{n-1}}. $$

\end{rem}

\noindent
We deduce a corollary:

\begin{corr} \label{cor}
Let $U$ a bounded \lip domain, for each $\epsilon >0$
there exists $V$ a $C^1$ domain such that,
$S = U \setminus \bar V \cup V \setminus \bar U$,
is a rectifiable open set verifying:
\begin{eqnarray} \label{approxopenset}
\mu (S) < \epsilon \\
\mathcal{H}^{n-1}(\partial S)<\epsilon.
\end{eqnarray}
\end{corr}

\noindent
\emph{Proof:}
We just present the main points for the proof of the corollary.
\\
By compacity of $\partial U$, there exists a finite open covering $(V_1,\ldots,V_k)$ of $\partial U$, 
such that for each open set we can trivialize the boundary. On each open set $V_i$,
we have by previous theorem a Lipschitz application $g_i: \partial U \cap V_i \mapsto \ms{R}^n$ which gives a $C^1$ 
hypersurface. We have $$\mathcal{H}^{n-1}(\{x \in \partial U \cap V_i | x \neq g(x)\}) \leq \epsilon.$$
\\
Moreover we can assume that this covering satisfies the following property.
Let $0<\eta<\epsilon$ and $Z:= \cup_{i \neq j} V_i \cap V_j$,
$\mathcal{H}^{n-1}(\partial U \cap  Z) \leq \eta$. 
We thus obtain an application $g: \partial U \mapsto \ms{R}^n$ which is Lipschitz
($\partial U$ is endowed with the induced metric by the euclidean metric on $\ms{R}^n$)
 and is the boundary of a $C^1$
domain $V$, for which we have:
$$\mathcal{H}^{n-1}(\{x \in \partial U | x \neq g(x)\}) \leq (k+1) \epsilon.$$
Then, $\mathcal{H}^{n-1}(\partial S)<2 \; \text{Lip}(g)^{n-1} (k+1) \epsilon$. 
And also with the same argument given in the preceding remark, there exists a constant $K$ such that,
$\mu(S) \leq K \; \text{Lip}(g) \left((k+1)\epsilon \right)^{1+\frac{1}{n-1}} $,
with $K=\sqrt{n-1} \; \text{Lip}(\partial U)$.
$\square$

\vspace{0.3cm}

\noindent
We now turn to the proof of the \textbf{lemma} \ref{lip1}.

\vspace{0.3cm}

\noindent
\emph{Proof:}
We use the corollary \ref{cor}, let $U_{\epsilon}$ a $C^1$ domain for $\epsilon$ as in the corollary. 
Let $M_1$ a constant such that in a compact neighborhood of $U$,
$|f \circ \phi_t - f | \leq M_1 t$, $g \leq M_2$, $|f| \leq M_3$ and $|X| \leq K$.
We have, with $$S_{\epsilon}=\Delta(U_{\epsilon},U)=U_{\epsilon} \setminus \bar U \cup U \setminus \bar U_{\epsilon},$$ 
We denote by $\theta = \mt{1}_{U}-\mt{1}_{U_{\epsilon}}$, so we have (triangular inequality for the second inequation):

\begin{eqnarray*}
|(J_t(U)-J_0(U))-(J_t(U_{\epsilon})-J_0(U_{\epsilon}))| &\leq & \int_{V} | \theta \circ \phi_t^{-1} f \circ \phi_t^{-1} g - \theta fg | d\mu, \\
|(J_t(U)-J_0(U))-(J_t(U_{\epsilon})-J_0(U_{\epsilon}))| &\leq & \int_{V} | \theta \circ \phi_t^{-1} (f \circ \phi_t^{-1} g-fg) | d\mu + \\
\int_{V} |(\theta \circ \phi_t^{-1}- \theta) fg | d\mu, \\
|(J_t(U)-J_0(U))-(J_t(U_{\epsilon})-J_0(U_{\epsilon}))| &\leq & \int_{V \cap \phi_t(S_{\epsilon})} |f-f \circ \phi_t^{-1}| |g|\mt{1}_V d\mu \\
+ \int_{V \cap \Delta(\phi_t(S_{\epsilon}),S_{\epsilon})} |f g| \mt{1}_V d\mu, \\
|(J_t(U)-J_0(U))-(J_t(U_{\epsilon})-J_0(U_{\epsilon}))| &\leq & t M_1 M_2 \mu(\phi_t(S_{\epsilon}) \cap V) + M_3 M_2 \mu(V \cap \Delta(\phi_t(S_{\epsilon}),S_{\epsilon}))
\end{eqnarray*}
We first treat the last term.
We claim that, for $s_0 > 0$ such that $\text{Lip}(\phi_t) \leq 2$, we have, for
$t \in [-s_0,s_0]$, $$\mu(\Delta(\phi_t(S_{\epsilon}),S_{\epsilon})) \leq t \max(2,M)^n \mathcal{H}^{n-1}(\partial S_{\epsilon})).$$  
Introduce $\Psi: (t,x) \in [-s_0,s_0] \times \ms{R}^n \mapsto \phi_t(x) \in \ms{R}^n$.
\\
We have $\text{Lip}(\Psi) \leq \max(2,M)$, and $\mathcal{H}^{n}([0,t] \times \partial S_{\epsilon})) = t \mathcal{H}^{n-1}(\partial S_{\epsilon})$.
Hence, $\mathcal{H}^{n}(\Psi([-s_0,s_0] \times \partial S_{\epsilon})) \leq t \max(2,M)^n \mathcal{H}^{n-1}(\partial S_{\epsilon}))$.
To finish, we prove that:
$$\Delta(\phi_t(S_{\epsilon}),S_{\epsilon}) \subset \bigcup_{s \leq t} \phi_t(\partial S_{\epsilon}).$$
Let $z \in \Delta(\phi_t(S_{\epsilon}),S_{\epsilon})$, 
\begin{itemize}
\item Suppose $z \notin S_{\epsilon}$, there exists $x \in S_{\epsilon}$ such that $\phi_t(x) = z$.
The map $c:s \in [0,t] \mapsto \phi_s(x)$ verifies $c(0)=x \in S$ and $c(t)=z \notin S$. By connexity,
there exists $u \in [0,t]$, such that $c(u) \in \partial S_{\epsilon}$. By composition of flow,
$\phi_{t-u}(c(u))=z$. 
\item Suppose $z \notin \phi_t(S_{\epsilon})$, there exists $x \in \phi_t(S_{\epsilon})$ such that $\phi_{-t}(x) = z$.
The map $c:s \in [0,t] \mapsto \phi_{-s}(x)$ verifies $c(0) \in \phi_t(S_{\epsilon})$ and $c(t)=z \notin \phi_t(S_{\epsilon})$.
By connexity, there exists $u \in [0,t]$, such that $c(u) \in \partial \phi_t(S_{\epsilon})$. By composition of flow,
$m=\phi_{-u}(z) \in \partial S_{\epsilon}$ and obviously, $\phi_{u}(m) = z$.  
\end{itemize}

\vspace{0.3cm}

\noindent
We give a bound for the first term in the same neighborhood for $t \in [-s_0,s_0]$,
$$\mu(\phi_t(S_{\epsilon}) \cap V) \leq \mu(\phi_t(S_{\epsilon})) \leq 2^n t \epsilon \mu(S_{\epsilon}).$$

\vspace{0.3cm}

\noindent
Consequently,
$$\lim \sup_{t\mapsto 0^+} | \frac{1}{t} [(J_t(U) -J_0(U) - (J_t(U_{\epsilon}) - J_0(U_{\epsilon}))] | \leq M_1 M_2 \mu(S_{\epsilon}) + M_3 \mathcal{H}^{n-1}(\partial S_{\epsilon})). $$
We can now obtain the conclusion.
Let $\epsilon>0$, 

\begin{equation*}
\lim  \sup_{t\mapsto 0^+} | \frac{1}{t} [(J_t(U) -J_0(U) - (J_t(U_{\epsilon}) - J_0(U_{\epsilon}))] | \leq  (M_1 M_2+M_3) \epsilon. \\
\end{equation*}

\noindent
We use now the formula already demonstrated for $C^1$ domains,
\begin{eqnarray*}
\lim  \sup_{t\mapsto 0^+} | \frac{1}{t} [(J_t(U) -J_0(U) - \int_{U_{\epsilon}} - <\nabla f,X> g \mt{1}_V d\mu + \\
\int_{\partial U_{\epsilon}} <X,n> fg  \tilde{\mt{1}}_V(X) d\mu_{|\partial U_{\epsilon}}] | \leq  (M_1 M_2+M_3) \epsilon, \\
\end{eqnarray*}
and the result is proven. 
$\square$

\section{Appendix}

\subsection{Central lemma of \cite{hamcurves}}

We present here a different version of the lemma, which is essentially the same, but from another point of view.
\begin{lm} \label{lmc}
Let $H$ a Hilbert space and $B$ a non-empty bounded subset of $E$ a Hilbert space such that
there exists a continuous linear application $g: H \mapsto E$. Assume that
for any $a \in H$, there exists $b_a \in B$ such that $\langle b_a, g(a) \rangle \geq 0$. 
Then, there exists $b \in \overline{\Conv(B)}$ such that $\langle b , g(a) \rangle = 0, \; \forall a \in H$.  
\end{lm}
\noindent
\emph{Proof:}
We denote by $H_0 = \overline{g(H)}$.
Let $p$ the orthogonal projection on $H_0$ and $Z=\overline{\Conv(B)}$ is a non-empty closed bounded convex subset of $H$,
whence weakly compact. As $p$ is weakly continuous and linear, $p(Z)$ is a weakly compact convex subset and thus strongly closed. 
From the projection theorem on closed convex subset, there exists $b \in p(Z)$ such that:
$|b| = \inf_{c \in Z} |c|$ and $\langle b,c-b \rangle \geq 0$ for $c \in p(Z)$. As a direct consequence,
we have also: 
\begin{equation} \label{ineqconv}
\langle b,u-b \rangle \geq 0 \; \forall u \in Z. 
\end{equation}
The element $b$ lies in the adherence of $g(H)$ then there exists a sequence $b_n \in g(H)$ such that 
$\lim b_n = b$. From the hypothesis, there 
exists $u_n \in B$ such that: $\langle u_n,-b_n \rangle \geq 0$.
As $B$ is bounded,
$\lim \langle u_n,b-b_n \rangle = 0$, hence $\limsup \langle u_n,b_n \rangle \leq 0$.
By \eqref{ineqconv}, we get $\limsup \langle u_n,b \rangle \geq \langle b,b \rangle$.
As a result, $\langle b,b \rangle \leq 0$ and $b=0$.
By definition, there exists $v \in Z$, $p(v)=b=0$. Now, $v \in H_0^{\perp}$ and $\phi(.)= \langle v,.\rangle$
gives the result.
$\square$

\subsection{A short lemma}

We give here a short proof of the perturbation of the flow of a time dependent vector field with respect to the vector field.
We assume in the proposition that the involved vector fields are $C^1$ but it can be proven with weaker assumptions on the regularity
of vector fields. (See \cite{JoanPhD}, for a detailed proof following another method.)

\begin{lm} \label{depder} 
Let $u_t$ and $v_t$ be two time dependent $C^1$ vector fields on $\ms{R}^n$, and denote by 
$\phi^{\epsilon}_{0,t}$ the flow generated by the vector field $u_t+\epsilon v_t$, then we have:
$$ \partial_{\epsilon} \phi_{0,1}(x) = \int_0^1 [d\phi_{t,1}]_{\phi_{0,t}(x)} v(\phi_{0,t}(x))dt.$$
\end{lm}

\noindent
\emph{Proof:}
Introduce the notation $A_t \in \ms{R}^n$ defined by:
$A_t(\phi_{0,t}(x)) = \partial_{\epsilon} \phi^{\epsilon}_{0,t}(x)$.
Deriving this expression with respect to the time variable:
$$ \frac{d}{dt} A_t(\phi_{0,t}(x)) = du_t(A_t(\phi_{0,t}(x))) + v_t(\phi_{0,t}(x)), $$
Remark that the expression above can be written as (with $\mathcal{L}$ the Lie derivative):
$$\mathcal{L}_{u_t} A_t = \frac{d}{du}_{|u=0} [d\phi_{0,t+u}]^{-1}(A_{t+u}(\phi_{t+u}(x)))=[d\phi_{0,t}]^{-1}(v_t(\phi_{0,t}(x))).$$
By integration in time, we obtain the result.
$\square$

\nocite{0915.49003}
\nocite{arn78}
\nocite{MR2062547}
\nocite{1062.93001}
\bibliographystyle{alpha}
\bibliography{biblio_art_1}

\end{document}